\documentclass[a4paper, 11pt, reqno]{amsart}
\usepackage{amsthm,amsmath,amssymb,mathrsfs,url,amsaddr}
\usepackage{datetime2}

\theoremstyle{definition}
\newtheorem{theorem}{Theorem}[section]
\newtheorem*{theorem*}{Theorem}
\newtheorem{lemma}[theorem]{Lemma}
\newtheorem{proposition}[theorem]{Proposition}
\newtheorem{corollary}[theorem]{Corollary}

\newtheorem{remark}[theorem]{Remark}

\newcommand\DN{\newcommand}
\DN\st{\,;\,}

\DN\lref[1]{Lemma~\ref{#1}}
\DN\tref[1]{Theorem~\ref{#1}}
\DN\pref[1]{Proposition~\ref{#1}}
\DN\sref[1]{Section~\ref{#1}}
\DN\dref[1]{Definition~\ref{#1}}
\DN\rref[1]{Remark~\ref{#1}} 
\DN\corref[1]{Corollary~\ref{#1}}
\numberwithin{equation}{section}
\newcounter{Const} 

\DN\Ct[1]{\refstepcounter{Const}c_{\theConst}\label{#1}}
\DN\cref[1]{c_{\ref{#1}}}

\newcommand{\ep}{\varepsilon}

\DN{\limi}[1]{\lim_{#1\to\infty}} 	

\DN\OD[2]{\frac{d #1 }{d #2}}
\DN\PD[2]{\frac{\partial #1 }{\partial #2}}
\DN\map[3]{#1:#2 \to #3}

\DN\R{\mathbb{R}}
\DN\N{\mathbb{N}}
\DN\Q{\mathbb{Q}}
\DN{\Z}{\mathbb{Z}}
\DN\C{\mathbb{C}}

\DN\www{\mathsf{w}}
\DN\xx{\mathbf{x}}
\DN\XX{\mathbf{X}}
\DN\xxx{\mathsf{x}}
\DN\yy{\mathbf{y}}
\DN\YY{\mathbf{Y}}
\DN\yyy{\mathsf{y}}
\DN\mm{\mathbf{m}}
\DN\zzz{\mathsf{z}}
\DN\pp{\mathbf{p}}
\DN\uu{\mathbf{u}}
\DN\uuu{\mathsf{u}}
\DN\eee{\mathsf{e}}
\DN\sss{\mathsf{s}}

\DN\ctil{\tilde{c}}

\DN\kk{\mathbf{k}}
\DN\E{\mathbb{E}}
\DN\MMM{\mathfrak{M}}
\DN\LLL{\mathfrak{L}}
\DN\LL{\mathscr{L}}

\DN\nWN{W_{\ge}^{N}}
\DN\nWNd{W_{\ge, d}^{N}}

\DN\PN{\mathscr{P}_N}
\DN\ad{\mathrm{ad}}
\DN\ip{\mathrm{ip}}
\DN\sym{\mathrm{sym}}

\DN\oFz{{}^{}_{1} F _{0} ^\theta }
\DN\oFo{{}^{}_{1} F _{1} ^\theta }
\DN\zFz{{}^{}_{0} F _{0} ^\theta }
\DN\zFo{{}^{}_{0} F _{1} ^\theta }
\DN\pFq{{}^{}_{p} F _{q} ^\theta }
\DN\one{\mathsf{1} }

\begin{document}

\title[]{Transforms between multivariate Laguerre and Meixner polynomials and their applications to gateways}

\author{KAWAMOTO Yosuke}
\address{Graduate school of environmental, life, natural science and technology, Okayama University, Okayama, Japan}
\email{y-kawamoto@okayama-u.ac.jp}

\author{SHIBUKAWA Genki}
\address{Kitami Institute of Technology, Hokkaido, Japan}
\email{g-shibukawa@mail.kitami-it.ac.jp}



\maketitle

\begin{abstract}
We establish gateways in both directions between the Laguerre process and a discrete Markov process for every inverse temperature $2\theta\ge 1$. 
Here, a gateway is an intertwining relation between Markov processes with continuous and discrete state spaces.
A key ingredient is a pair of generating function identities that relate two double generating functions for multivariate Meixner polynomials to a single and a double generating function for multivariate Laguerre polynomials, respectively.
These identities reveal a close relationship between the two families of polynomials and yield integral transforms that underlie the gateway relations.
Gateway relations provide a means of transferring properties between continuous and discrete processes. 
As an application, we derive the Feller--Dynkin property and the uniqueness of the invariant distribution for the Laguerre process from the corresponding properties of the discrete process.
\end{abstract}

\tableofcontents

\section{Introduction}
\subsection{Gateway between non-colliding processes}
An intertwining relation means that two stochastic processes commute through a Markov kernel linking their state spaces.
More precisely, given two Markov semigroups $\{T_t^A \}_{t\ge 0 }$ and $\{T_t^B \}_{t\ge 0 }$ on state spaces $A$ and $B$, respectively, and a Markov kernel $\Lambda : A \dashrightarrow B$, an intertwining relation is the identity  
\begin{align*}
 T_t^A \Lambda =\Lambda T_t^B 
\end{align*}
for all $t \ge 0$.
Intertwining relations have proved useful in the analysis of stochastic processes. 
Early important examples appear in the work of Rogers and Pitman \cite{RoP81} on Markov functions, while Diaconis and Fill \cite{DiF90} used intertwining to investigate the convergence of Markov chains to stationarity.

Intertwining relations have been extensively studied for stochastic processes arising in random matrix theory.
A typical example is a system of non-intersecting Brownian motions, also known as the Dyson Brownian motion at inverse temperature $2$ \cite{War07}.
Intertwining relations for the Laguerre process at inverse temperature $2$, an analogous non-intersecting diffusion, were obtained in \cite{AOW19, BuK25}.
For general inverse temperature, the intertwining relations for the Dyson Brownian motions were established in \cite{GoS15,RaS18}, and those for the Laguerre processes were obtained in \cite{Ass19a, KaS26}.
In parallel, intertwining relations for continuous-time Markov chains related to random matrix theory have also been investigated, notably for multilevel dynamics on Gelfand--Tsetlin patterns and related interlacing particle systems \cite{Ass18, BoO12, BoO13, WaW09}.

An intertwining relation that directly connects a continuous Markov process with a discrete one, as an alternative to approximation procedures, is called a gateway in \cite{MiP19}, and we adopt this terminology.
In the context of random matrix theory, gateways have been less extensively studied than intertwining relations between diffusions or between discrete Markov chains.
In particular, whilst intertwining relations are known for diffusion processes at general inverse temperature, we are not aware of a gateway in this generality.
In this paper, using an approach based on symmetric polynomials, we establish gateways for Laguerre processes and the corresponding Markov chains at general inverse temperature.

We set $\R_{\ge} =[0,\infty) $ and $\Z_{\ge} =\N \cup \{0\}$.
The closed non-negative Weyl chamber and the discrete Weyl chamber are defined by 
\begin{align*}
\nWN &=\{ \xxx =(x_1,\ldots,x_N) \in \R_{\ge}^N \st  x_1\le x_2 \le \cdots \le x_N  \}
,  \\
 \nWNd &=\{\xx =(x_1,\ldots,x_N) \in \Z_{\ge}^N \st x_1< x_2 < \cdots <x_N \}  
.\end{align*}
To distinguish continuous objects from discrete ones, we use sans-serif symbols for the former and boldface symbols for the latter.
Let $\Delta(\xxx )=\prod_{1\le i<j\le N}(x_j-x_i)$ for $\xxx \in \R^N$.

We begin by recalling a known gateway between two non-colliding processes.
For $\alpha >0$, consider the diffusion  on $\mathbb{R}_{\ge }$ generated by the differential operator
\begin{align*}
\mathsf{A}_{\alpha} =x\frac{d^2}{dx^2 }+ (-x+\alpha )\frac{d}{dx} 
.\end{align*}
Let $\mathsf{p}_{\alpha,t}(x,y)$ denote its transition density.
Consider the non-colliding system of $N$ independent copies of the diffusion generated by $\mathsf{A}_\alpha$. 
Its state space is $\nWN$, and its transition density is given by the Karlin--McGregor determinantal formula
\begin{align*}
\mathsf{p}_{\alpha , t}^{N}(\xxx ,\yyy )
=e^{\frac{N(N-1)}{2}t}\frac{\Delta (\yyy ) }{\Delta (\xxx ) }
\det [\mathsf{p}_{\alpha, t}(x_i,y_j)]_{i,j=1}^N
\end{align*}
for $t>0$ and $x_1<\cdots<x_N$.
Let $\{ \mathsf{T}_{\alpha, t}^N\}_{t\ge 0}$ be the Markov semigroup associated with this non-colliding diffusion.

We next introduce the corresponding discrete process.
Let $\nabla_{+}$ and $\nabla_{-}$ denote forward and backward difference operators,  defined by $\nabla_{\pm }f (x)=f(x\pm1 )-f(x)$.
For $\alpha, \sigma >0$, define the difference operator 
\begin{align*}
  \mathbf{A}_{\alpha, \sigma} =\sigma(x + \alpha  )\nabla_{+} +(\sigma+1)x \nabla_{-}
,\end{align*}
which is the generator of a continuous-time Markov chain on $\Z_{\ge }$.
Let $\mathbf{p}_{\alpha, \sigma , t}^{} (x,y)$ be the transition probability of the Markov chain generated by $\mathbf{A}_{\alpha, \sigma}$.
Then, consider the non-colliding system of $N$ independent copies of this chain.
Its state space is $\nWNd$, and its transition probability is also given by the Karlin--McGregor formula
\begin{align*}
\mathbf{p}_{\alpha,\sigma , t}^{N}(\xx,\yy)=e^{\frac{N(N-1)}{2}t}  \frac{\Delta(\yy) } {\Delta (\xx ) } \det [\mathbf{p}_{\alpha, \sigma, t}^{} (x_i,y_j)]_{i,j=1}^N
.\end{align*} 
Let $\{ \mathbf{T}_{\alpha,\sigma , t}^{N} \}_{t\ge 0}$ be the Markov semigroup of this non-colliding chain.

We introduce a kernel linking the continuous space $\nWN $ and the discrete space $\nWNd $.
For $\sigma>0 $, define a Markov kernel $\Lambda_{\sigma}^N : \nWN \dashrightarrow \nWNd $ by 
\begin{align*}
\Lambda_{\sigma}^N (\xxx, \yy ) &= \frac{\Delta (\yy) } {\Delta(\sigma \xxx ) } \det\bigg[\frac{(\sigma x_i)^{y_j} e^{-\sigma x_i}}{y_j !}\bigg]_{i,j=1}^{N}
.\end{align*}
Strictly speaking, this expression is well-defined for $\xxx$ in the interior of $\nWN$ and admits a continuous extension to all $\xxx\in\nWN$, as discussed in \cite{Ass19b}.
Then, the following gateway relation was established in \cite[Theorem 1.8]{Ass19b}: assume that $\alpha, \sigma >0 $.
Then, for any $N\in \N$ and $t\ge 0$, we have the equality of Markov kernels
\begin{align}\label{:11b}
\mathsf{T}_{ \alpha ,t}^{N} \Lambda_{\sigma}^N
&
=\Lambda_{\sigma}^N\mathbf{T}_{\alpha,\sigma ,t}^N
.\end{align}
The relation \eqref{:11b} was proved for $N=1$ in \cite[Proposition 6.7]{BoO13} and \cite[Proposition 4.4]{MiP19}, and was then extended to general $N$ using the Karlin--McGregor formula.

\subsection{Gateways for general inverse temperature}

Let $\theta \ge 1/2 $ and $\alpha >0$.
The $N$-dimensional Laguerre process with parameter $\alpha$ and inverse temperature $2\theta$ is defined as a solution to the stochastic differential equation
\begin{align}\label{:11a}
  dX_t^{i} =\sqrt{2X_t^i } dB_t^i -X_t^idt + \theta \Big (\alpha + \sum_{j\neq i}^N \frac{2X_t^i}{X_t^i-X_t^j} \Big) dt
,\end{align}
for $i=1,\ldots,N$, where $B^1,\ldots,B^N$ are independent
standard Brownian motions.
For any starting point $\xxx \in \nWN$, equation~\eqref{:11a} has a unique strong solution with no collisions and no explosions by a  slight modification of the proof of \cite[Corollary 6.5]{GrM14}.
Let $\{ \mathsf{T}_{\theta, \alpha, t}^N \}_{t\ge 0}$ denote its Markov semigroup.

For $\theta=1$, the Laguerre process coincides with the non-colliding $\mathsf{A}_\alpha $-diffusions; see \cite{KoO01} for the case without the Ornstein--Uhlenbeck term.
It is therefore natural to consider extending the gateway \eqref{:11b} to general $\theta$, and one of the main purposes of this paper is to establish such an extension.

However, the proof of \eqref{:11b} relies on the Karlin--McGregor formulae, which are specific to $\theta=1$, and therefore the argument does not apply to $\theta$-extensions.
More fundamentally, we need to identify appropriate extensions of the kernel $\Lambda_{\sigma }^N$ and the discrete dynamics $\mathbf{T}_{\alpha, \sigma,t}^N$ in the parameter $\theta $.

For this purpose, an approach based on symmetric polynomials proves useful.
We use multivariate Laguerre and Meixner polynomials, and it turns out that certain identities for their generating functions underlie the $\theta$-extension of \eqref{:11b}.

\subsection{Meixner and Laguerre polynomials in one variable}\label{s:13}
Before working with symmetric polynomials, we illustrate our strategy in the one-variable case.
Let $(x)_{z}=\Gamma(x+z)/\Gamma(x)$ be the shifted factorial.
The Laguerre and Meixner polynomials are defined by
\begin{align*} 
  L_m^a(x)&
  :=\frac{(a +1)_m}{m!}\sum_{k=0}^m \frac{(-1)^k}{(a+1)_k} \binom{m}{k} x^k
,  \\
  M_{m} (x; a,c ) &:= \sum_{k=0}^{m} \frac{k!}{(a)_{k}} \binom{m}{k} \binom{x}{k} \Big(1-\frac{1}{c} \Big)^k 
.\end{align*}
Note that $L_m^{\alpha -1}(x)$ and  $M_{m} (x; \alpha , \frac{\sigma}{1+\sigma}  )$ are eigenfunctions of the generators $\mathsf{A}_{\alpha }$ and $\mathbf{A}_{\alpha, \sigma}$, respectively.
The key identity underlying the gateway \eqref{:11b} for $N=1$ is the following transform, proved by a direct calculation in \cite[(6.17)]{BoO13}:
\begin{align}  \label{:21a}
\sum_{y \in \Z_{\ge }} \frac{(\sigma x)^y e^{-\sigma x} }{y!}  M_{m} \Big(y; \alpha , \frac{\sigma}{1+\sigma} \Big) &= \frac{m!}{(\alpha )_{m }} L_{m}^{\alpha-1} (x)
.\end{align}

Here, we present an alternative proof of \eqref{:21a} using certain generating functions.
Equation \eqref{:21a} follows by comparing the coefficients of $u^m$ in the identity
\begin{align}  \label{:21f}
  \sum_{m \in \Z_{\ge} } \frac{1}{(\alpha )_m} L_m^{\alpha -1}(z) u^m
  &
  =e^u {_{0}F_1} \Big(\begin{matrix}  \\ \alpha  \end{matrix};-zu  \Big) 
  \\&\notag
  =\sum_{m,x \in \Z_{\ge} } M_m(x ; \alpha ,c)e^{- \frac{c}{1-c}z}  \frac{u^m}{m!} \frac{(\frac{c}{1-c}z)^x}{x!}
.\end{align}
The first equality is a generating function of the Laguerre polynomials \cite[(9.12.11)]{KLS10}, and the second one follows from applying the binomial formula. 
The last expression of \eqref{:21f} is a generating function of a generating function.
In this sense, we call it a double generating function for the Meixner polynomials.

Thus, we obtain a key identity: the generating function of the Laguerre polynomials coincides with the double generating function of the Meixner polynomials.
This identity reveals a close connection between these two families of polynomials and offers a useful perspective on their structure.
We obtain a multivariate analogue of \eqref{:21f}, which yields basic formulae for the Meixner polynomials such as a multivariable version of \eqref{:21a}.
As a result, we establish a $\theta $-extension of \eqref{:11b}.

The gateway considered above is realised by a Markov kernel from $\nWN $ to $\nWNd$.
There is also a gateway in the opposite direction.
We consider the identity
\begin{align}\label{:14a}
  \frac{1}{\Gamma(\alpha ) } \! \int_0^\infty \! \!e^u {_{0}F_1} \Big(\begin{matrix}  \\ \alpha \end{matrix};(c-1)y u \Big)  e^{-y(1-z) } y^{\alpha -1}dy
  &
  =\frac{e^{cu}}{(1-z)^{\alpha }}  \exp \Big( (c-1) \frac{uz}{1-z}\Big)
  \\& \notag 
=\sum_{m,x \in \Z_{\ge}} \! \!(\alpha )_x  M_{m} (x; \alpha ,c ) \frac{(cu)^m}{m! } \frac{z^x}{x!}
.\end{align}
The first expression is the (modified) Laplace transform of the Laguerre generating function \eqref{:21f}; it is reasonable to call it a double generating function for the Laguerre polynomials.
On the other hand, the last expression is a double generating function for the Meixner polynomials, different from the one in \eqref{:21f}.
Thus, the identity \eqref{:14a} shows that these two double generating functions coincide.
Comparing the coefficients of $u^mz^x$ in \eqref{:14a}, we obtain the integral transform 
\begin{align}
\label{:21e}
\frac{1 }{ \Gamma(\alpha +x) }\int_0^\infty L_{m}^{\alpha -1}((1-c)y ) e^{-y} y^{\alpha +x-1}  dy &= c^m \frac{ (\alpha )_m}{m!}  M_{m}(x; \alpha ,c) 
.\end{align}
This transform gives rise to an appropriate Markov kernel from $\Z_{\ge} $ to $\R_{\ge } $ for which the diffusion and the Markov chain under consideration satisfy a gateway.

\subsection{Main results}\label{s:14}
We extend the argument in \sref{s:13} to the multivariate case.
For this purpose, we use the multivariate Laguerre polynomials introduced in \cite{Las91}, whose basic properties were studied in \cite{BaF97}.
This family provides the appropriate multivariate version for the above argument.
In particular, the multivariate Laguerre polynomials are eigenfunctions of the generator of the Laguerre process.

The multivariate Meixner polynomials that are appropriate for our purposes were introduced and studied in \cite{Shi16} for $\theta=1/2,1, 2$ in the context of harmonic analysis on symmetric cones.
We investigate this family in detail for all $\theta>0$.
To this end, the multivariate versions of \eqref{:21f} and \eqref{:14a} play central roles.
These identities shed light on the relationship between multivariate Laguerre and Meixner polynomials, and they can be viewed as master formulae in this paper.
For example, we derive a difference operator for which Meixner polynomials are eigenfunctions.
The Markov chain generated by this operator, defined as follows, is in a gateway relation with the Laguerre process.

For $i=1,\ldots,N$, let $\eee_i$ denote the vector with $1$ in the $i$-th coordinate and $0$ in all other coordinates.
Let $ \nabla_{+ , i}$ and $ \nabla_{- , i}$ be the forward and backward difference operators defined by $(  \nabla_{\pm , i} f ) (\xxx ) =f(\xxx \pm \eee_i) - f(\xxx )$.
Let $\mathbf{A}_{\theta,\alpha ,\sigma }^{N}$ be the difference operator 
\begin{align}\label{:11g}
  \mathbf{A}_{\theta,\alpha ,\sigma }^{N}=\sum_{i=1}^N  \Big( \sigma q_{\theta, \alpha }^{N} (\xx, \xx+ \eee _i ) \nabla_{+,i} + (1+\sigma) q_{\theta, \alpha }^{N} (\xx, \xx- \eee_i ) \nabla_{-,i} \Big)
,\end{align}
where 
\begin{align}\notag
  q_{\theta, \alpha }^{N} (\xx, \xx+ \eee _i )&= (x_i-i+1+\theta (\alpha+i-1) ) \prod_{1\le j \neq i \le N } \frac{x_i-x_j-i+j+\theta(i-j+1) }{x_i-x_j-i+j+\theta(i-j) } 
,\\ \notag
  q_{\theta, \alpha }^{N} (\xx, \xx- \eee_i )&= (x_i-i+1+\theta(i-1) )\prod_{1\le j \neq i \le N } \frac{x_i-x_j-i+j+\theta(i-j-1)}{x_i-x_j-i+j+\theta(i-j) } 
.\end{align}
Then, the generator $\mathbf{A}_{\theta,\alpha ,\sigma }^{N}$ determines a continuous-time Markov chain on $\nWNd $; see \lref{l:51}.
For future reference, we call this chain the Meixner process.
Let $ \{\mathbf{T}_{\theta,\alpha, \sigma ,t}^{N}\}_{t\ge 0}$ be the associated Markov semigroup.

Motivated by the transforms corresponding to the multivariate analogues of \eqref{:21a} and \eqref{:21e}, given in \eqref{:27a} and \eqref{:36d}, respectively, we define Markov kernels $ \Lambda_{\theta, \sigma }^{N} : \nWN \dashrightarrow \nWNd $ and $\Lambda_{\theta, \alpha , \sigma}^{*, N} : \nWNd \dashrightarrow \nWN $, which will be specified in \eqref{:41c} and \eqref{:41d} by using the Jack polynomials.
Then we have the gateway relations between the Laguerre and Meixner processes.
Note that the following gateway \eqref{:12a} for $\theta=1$ coincides with \eqref{:11b}.

\begin{theorem}\label{t:11}
Assume that $\theta\ge 1/2$ and $\alpha, \sigma>0$.
Then, for any $N\in \N$ and $t\ge 0$, we have the equalities of Markov kernels
\begin{align}\label{:12a}
\mathsf{T}_{\theta, \alpha ,t}^{N} \Lambda_{\theta,\sigma}^N
&
=\Lambda_{\theta,\sigma}^N \mathbf{T}_{\theta, \alpha ,\sigma ,t}^N
,\\
\mathbf{T}_{\theta, \alpha , \sigma ,t}^N \Lambda_{\theta, \alpha ,\sigma}^{*,N}
 &\label{:12b}
=\Lambda_{\theta, \alpha ,\sigma}^{*,N} \mathsf{T}_{\theta,\alpha ,t}^{N}
.\end{align}
\end{theorem}

Analogues of \eqref{:12a} and \eqref{:12b} can also be obtained between the non-stationary versions of the Laguerre and Meixner processes.
These provide $\theta$-extensions of the gateways in \cite[Theorem 1.4]{Ass19b}. 
Their derivation requires further analysis of Meixner polynomials, and will be presented in a sequel to this paper.

Intertwining relations constitute a useful integrable structure and have found a variety of fruitful applications. 
In particular, gateway relations provide an exact connection between continuous and discrete processes, allowing one to transfer information between them. 
Specifically, discrete self-similarity and a discrete analogue of the beta-gamma algebra were obtained by transferring the  corresponding properties of a diffusion to a Markov chain through a gateway relation \cite{MiP19}. 
A related gateway idea was later developed to derive scaling limits and ergodic properties for discrete self-similar and ergodic Markov chains \cite{MPS22}.
In this paper, as applications of the gateways, we derive the uniqueness of the invariant distribution and the Feller--Dynkin property for the Laguerre processes from the corresponding properties of the discrete processes.

\subsection{Organisation of the paper}
\sref{s:2} introduces Jack and multivariate Laguerre polynomials.
In \sref{s:3}, we establish basic formulae for multivariate Meixner polynomials.
\sref{s:4} collects some facts about Markov kernels and continuous-time Markov chains. 
\sref{s:5} is devoted to proving \tref{t:11}.
In \sref{s:6}, we express our general results for $\theta=1$ in terms of well-known Schur-type multivariate polynomials.
In Appendix~\ref{s:A}, we collect formulae for multivariate Charlier and Krawtchouk polynomials, which are obtained from that for Meixner.
In Appendix~\ref{s:B}, we collect estimates for special functions.

\section{Preliminaries on Jack and multivariate Laguerre polynomials}\label{s:2}
Let $\PN=\{\mm=(m_1, \dots , m_{N})\in \Z_{\ge}; m_1 \ge \dots \ge m_N\ge0  \} $.
For $\mm \in \PN$, we set $\mm !=m_1 ! \cdots m_N!$.
For $\kk, \mm \in \PN $, we write $\kk \subseteq \mm $ if $k_i \le m_i$ holds for all $i$.
Let $\one $ denote the $N$-dimensional vector whose entries are all $1$, and we write $\delta=(N-1, N-2,\ldots, 1,0 ) $.
For $x\in\R$ and $\mm\in\PN$, we write $(x)_\mm^\theta  =\prod_{j=1}^{N} (x-\theta(j-1))_{m_j}$.
For $\xxx =(x_1,\ldots, x_N)\in \R^N  $, we write $|\xxx | =\sum_{i=1}^{N} x_i$ and define the supremum norm $\|\xxx \|=\sup_{1\le i\le N} |x_i |$.
Let $\mathfrak{S}_N$ be the symmetric group of degree $N$.
We denote the ring of symmetric polynomials over a field $K$ by 
\begin{align*}
  K [\xxx ]^{\mathfrak{S}_N}=\{f\in  K [\xxx ] \st f(x_{\sigma(1)},\ldots , x_{\sigma(N)}) =f(x_1,\ldots , x_N) \text{ for all } \sigma \in \mathfrak{S}_N \}
.\end{align*}

Hereafter, unless otherwise stated, we assume that $a \in \R$ satisfies $(a)_\mm^\theta \neq 0$ for any $\mm \in \PN$.
Furthermore, throughout this section and \sref{s:3}, we assume $\theta >0$.

\subsection{Jack polynomials}
We begin by introducing Jack polynomials.
See \cite[Chapter 2]{ViK95} for further details.
We introduce the following differential operators in $N$-variables
\begin{align*}
E_{k} ^N := E_{k } ^{N,(\xxx )}&:=\sum_{i=1}^{N} x_i^k \PD{}{x_i}, 
\\
 D_{k} ^N:=D_{k} ^{N,(\xxx )} &:= \sum_{i=1}^{N} x_i^k \frac{\partial ^2 }{ \partial x_i^2 } + 2\theta \sum_{i=1}^{N} \sum_{1\le j \neq i\le N} \frac{x_i^k}{x_i-x_j} \PD{}{x_i}
.\end{align*}

 For $\mm  \in\PN $, the Jack polynomial $ P_\mm ^\theta $ is the unique symmetric polynomial satisfying $$ D_2^N  P_\mm ^\theta (\xxx ) = \Big( \sum_{i=1}^N m_i(m_i-1) +2\theta \sum_{i=1}^N(N-i)m_i \Big) P_\mm ^\theta (\xxx )$$
and having the form  $ P_\mm ^\theta  (\xxx ) =m_\mm+\sum_{\kk<\mm}c_{\mm\kk}m_\kk $, where $<$ denotes the dominance order and $m_\mm$ is the monomial symmetric polynomial.
We shall use the following formula; see \cite[(4.8)]{Koo15} and \cite[(6.4)]{OkO97}:
\begin{align}\notag
P_\mm ^\theta(\one  ) &= \prod_{1\le i <j \le N}(m_i-m_j +\theta(j-i))_\theta  \prod_{k=1}^{N} \frac{\Gamma(\theta)}{\Gamma(\theta k)}
\\ \label{:20a}
&=\prod_{1\le i<j\le N} \frac{(\theta(j-i+1))_{m_i-m_j}}{(\theta(j-i))_{m_i-m_j}}  
.\end{align}

The tableau sum formula for the Jack polynomials \eqref{B:1c} immediately yields the following non-negativity result.
\begin{lemma}\label{l:1}
For $\theta>0$, $\mm\in\N$, and $\xxx \in \nWN $, we have  $P_{\mm }^\theta(\xxx ) \ge 0$.
\end{lemma}

We next introduce the interpolation Jack polynomial $P_\mm^{\theta,\ip}$; see \cite{KnS96, Koo15, Oko98, OkO97, Sah94}.
It is the unique $\mm$-monic symmetric polynomial of degree $|\mm |$ such that $P_\mm^{\theta,\ip } (\kk + \theta \delta )=0$ for each partition $\kk \neq \mm $ with $|\kk | \le |\mm |$.
Note that the highest-degree homogeneous component of $P_\mm^{\theta,\ip }$ is $P_{\mm} ^{\theta}$, and that we have
\begin{align}\label{:20b}
P_\mm^{\theta,\ip } (\mm + \theta\delta ) =\prod_{k=1}^N ( \theta(N-k)+1 )_{m_k} \prod_{1\le i<j\le N} \frac{(\theta(j-i-1)+1)_{m_i-m_j}}{(\theta(j-i)+1)_{m_i-m_j}}  
,\end{align}
see \cite[(7.4), (7.5)]{Koo15}.
Define two Jack polynomials with different normalisations by 
\begin{align}\label{:20c}
 \Phi _{\mm }^{\theta }(\xxx )
 =\frac{P_{\mm }^\theta (\xxx )}{P_{\mm }^\theta (\one  )}, 
 \qquad 
  \Psi _{\mm }^{\theta }(\xxx)=\frac{P_{\mm }^\theta (\xxx )}{P_{\mm }^{\theta,\ip } (\mm +\theta \delta)} 
.\end{align}

For $\xxx \in \R^N$ and $\kk \in \PN $, we introduce the Jack binomial coefficient by
\begin{align*}
 \binom{\xxx }{\kk }_{\theta } = \frac{P_\kk^{\theta, \ip } (\xxx + \theta \delta )}{P_\kk ^{\theta, \ip} (\kk + \theta \delta  ) } 
.\end{align*}
For  all $\mm,\kk \in \PN $, we have $ \binom{\mm }{\kk }_\theta  \ge 0  $ \cite{Sah11}, and, in particular,  $\binom{\mm }{\kk }_{\theta }=0$ unless $\kk \subseteq \mm $.
The Jack binomial formula was established in \cite{Las90}:
\begin{align}\label{:31f}
\Phi_{\mm }^\theta (\one +\xxx )=\sum_{\kk \subseteq \mm }  \binom{\mm }{\kk }_{\theta }\Phi_{\kk } ^\theta (\xxx )
.\end{align}

Let $a_1,\ldots, a_p, b_1,\ldots, b_q \in \R $ be such that none of the numbers $b_j-\theta (i-1)$ ($j=1,\ldots, q, i=1,\ldots, N$) coincides with a nonpositive integer.
Define multivariate hypergeometric functions associated with Jack polynomials (see 
\cite[Chapter 3]{ViK95}):
\begin{align}\label{:21v}
 \pFq \Big(\begin{matrix} a_1 , \ldots, a_p  \\ b_1 , \ldots , b_q  \end{matrix} \st \xxx, \yyy \Big)=
  \sum_{\mm \in \PN }  \frac{(a_{1} )_{\mm }^\theta \cdots (a_p )_{\mm }^\theta  }{(b_1 )_{\mm } ^\theta \cdots (b_{q} )_{\mm }^\theta }\Phi_{\mm }^\theta(\xxx ) \Psi_{\mm}^\theta(\yyy )
.\end{align}

\begin{lemma}\label{l:B2}\cite[Theorem 6.5, Remark 6.7]{BrR23}
  
\begin{enumerate}
\item[(i)]
The hypergeometric function \eqref{:21v} converges absolutely for all $\xxx, \yyy \in\R^N $ if $p\le q$, and for $\xxx, \yyy $ with $\| \xxx \| \| \yyy \|<1 $ if $p=q+1$.

\item[(ii)]
Assume that  $p\le q$ and $\ep >1$.
Then, there exits a constant $C_\ep >0$ such that,
\begin{align}\label{B:2c}
\sum_{\mm \in \PN } \bigg| \frac{(a_{1} )_{\mm }^\theta \cdots (a_p )_{\mm }^\theta  }{(b_1 )_{\mm } ^\theta \cdots (b_{q} )_{\mm }^\theta }\Phi_{\mm }^\theta(\xxx ) \Psi_{\mm}^\theta(\yyy ) \bigg|
\le  e ^{ N \ep \|\xxx \| \|\yyy \| } 
.\end{align}  
\end{enumerate}

\end{lemma}

\subsection{Multivariate Laguerre polynomials}
Here we collect some facts about multivariate Laguerre polynomials; see \cite{BaF97} for further details.
Define multivariate Laguerre polynomials by
\begin{align}\label{:22a}
  \LL_{\mm}^{\theta,a}(\xxx )& = \Psi_{\mm}^{\theta}( \one ) \sum_{\kk \subseteq \mm }  (-1)^{|\kk |}
  \binom{\mm }{\kk }_{\theta } \frac{(a )_{\mm }^\theta }{(a  )_{\kk} ^\theta} \Phi_{\kk }^{\theta }(\xxx )
.\end{align}
We remark that our definition is different from that in \cite{BaF97}. 
The multivariate Laguerre polynomial used in \cite{BaF97}, denoted by $L_{\mm }^{a }(\xxx ; \theta^{-1} )$, is related to ours by
\begin{align*}
L_{\mm }^{a }(\xxx ; \theta^{-1} )=\frac{1 }{|\mm |! \Psi_{\mm}^{\theta}( \one  ) } \LL_{\mm }^{\theta , a+1+\theta(N-1)} (\xxx ) 
.\end{align*}
In particular, for $N=1$, we have $\LL_{m }^{\theta, a } (x )=L_{m}^{a-1} (x )$.

Introduce the second differential operator
\begin{align*}
\mathsf{A}_{\theta, a} ^N := \mathsf{A}_{\theta, a } ^{N,(\xxx )} := D_{1}^{N,(\xxx )} - E_{1}^{N,(\xxx )} + (a-\theta(N-1)) E_{0}^{N,(\xxx )}
.\end{align*}
Then, Laguerre polynomials satisfy
\begin{align}\label{:31a}
 \mathsf{A}_{\theta, a} ^N \LL_{\mm}^{\theta,a}(\xxx ) &= -|\mm | \LL_{\mm}^{\theta,a}(\xxx ) 
.\end{align}

Multivariate Laguerre polynomials admit three types of generating functions.
In the following lemma, for $\uuu \in\R^N$, we interpret 
\begin{align*}
\frac{\uuu }{\one - \uuu }=\Big(\frac{u_1 }{1- u_1 },\ldots ,\frac{u_N }{1- u_N } \Big)  
,\end{align*}
and use similar componentwise notation throughout.
\begin{lemma}\cite[Proposition 4.1, Proposition 4.2]{BaF97} \label{l:21}
For any $\gamma \in \R $, we have 
\begin{align}
\label{:26d}
 \sum_{\mm \in \PN } \frac{(\gamma )_{\mm} ^\theta }{(a)_{\mm} ^\theta }\LL_{\mm}^{\theta,a}(\zzz ) \Phi_{\mm }^\theta (\uuu )
=\prod_{j=1}^{N} (1-u_j)^{-\gamma } \oFo \Big(\begin{matrix} \gamma \\ a  \end{matrix};-\zzz , \frac{\uuu }{\one - \uuu } \Big) 
.\end{align}
In particular, we have 
\begin{align} \label{:26c}
\sum_{\mm \in \PN } \LL_{\mm}^{\theta,a}(\zzz ) \Phi_{\mm }^\theta (\uuu )
  &=\prod_{j=1}^{N} (1-u_j)^{-a} \zFz \Big(;-\zzz, \frac{\uuu }{\one -\uuu }  \Big) 
 ,\\ \label{:26a}
\sum_{\mm \in \PN } \frac{1}{(a)_{\mm } ^\theta} \LL_{\mm }^{\theta,a}(\zzz ) \Phi_{\mm }^\theta (\uuu )
  &=e^{|\uuu |} \zFo \Big(\begin{matrix}  \\ a  \end{matrix};-\zzz , \uuu  \Big) 
.\end{align} 
\end{lemma}

We remark that \eqref{:26d} and \eqref{:26c} holds for $\uuu \in \R^{N}$ with sufficiently small $\|\uuu \|$ and all $\zzz \in \R^N$, whereas \eqref{:26a} holds for all $\uuu , \zzz \in \R^N$.

Clearly, equation \eqref{:26c} is a special case of \eqref{:26d}.
Furthermore, equation \eqref{:26a} can be obtained by replacing $\uuu $ by $\uuu /\gamma $ in \eqref{:26d} and taking $\gamma \to \infty$ with the degeneration  
\begin{align*}
\limi{\gamma }  \oFo \Big(\begin{matrix} \gamma  \\ a  \end{matrix} \st -\zzz , \frac{\gamma^{-1} \uuu  }{\one  -\gamma ^{-1} \uuu  }   \Big)   = \zFo  \Big(\begin{matrix}  \\ a  \end{matrix} \st -\zzz , \uuu  \Big) 
.\end{align*}

We note that the double series obtained by expanding the Laguerre polynomials in the series on the left-hand side of \eqref{:26d} converges absolutely and  locally uniformly in $\zzz$ for sufficiently small $\|\uuu\| $ by \lref{l:B3}.
Hence, rearranging the double series and comparing the coefficients yields the following summation formulae.
\begin{lemma}
For any $\kk\in \PN $, we have
\begin{align} \label{:28a}
\prod_{j=1}^{N} (1-u_j)^{-a} \Psi_{\kk }^\theta \Big(\frac{\uuu }{\one -\uuu } \Big)
& =
\sum_{ \mm \in \PN } \binom{\mm }{\kk }_{\theta } \frac{(a  )_{\mm }^\theta }{(a  )_{\kk}^\theta }
\Psi_{\mm }^\theta (\uuu )
, \\
\label{:25b}
  e^{|\uuu |} \Psi_{\kk }^\theta  (\uuu )
  &=\sum_{\mm \in \PN } \binom{\mm }{\kk }_\theta \Psi_{\mm }^\theta (\uuu )
.\end{align}
Here, \eqref{:28a} holds for $\uuu \in \R^{N}$ with sufficiently small $\|\uuu \|$, whereas \eqref{:25b} holds for any $\uuu \in \R^{N}$.
\begin{proof}
Equation \eqref{:28a} follows from \eqref{:26d} or \eqref{:26c}.
Indeed, expanding both sides of \eqref{:26c} by \eqref{:21v} and \eqref{:22a}, we obtain
\begin{align*}
\sum_{\mm \in \PN } \sum_{\kk \subseteq \mm }  (-1)^{|\kk |}
  \binom{\mm }{\kk }_{\theta } \frac{(a )_{\mm }^\theta }{(a  )_{\kk}^\theta } \Phi_{\kk }^{\theta}(\zzz )
\Psi_{\mm }^\theta (\uuu )
   =\sum_{\kk\in\PN } \prod_{j=1}^{N} (1-u_j)^{-a} \Phi_{\kk }^\theta (-\zzz ) \Psi_{\kk }^\theta \Big(\frac{\uuu }{\one -\uuu } \Big)  
\end{align*}
Comparing the coefficients of $\Phi_{\kk }^{\theta }(\zzz )$, we obtain \eqref{:28a}. 
Equation \eqref{:25b} readily follows from the Jack binomial formula \eqref{:31f}, or alternatively from a scaling limit of \eqref{:28a}.
\end{proof}
\end{lemma}

Hereafter, we often omit explicit conditions on $\uuu, \zzz $, and identities involving generating functions are understood to hold whenever both sides are well defined.
We state explicitly that an identity holds for all $\uuu\in\R^N$ only when this is needed in the arguments below.

\begin{lemma}
We have
\begin{align}\notag
 \oFz \Big(  \begin{matrix} a \\ \, \end{matrix} ;\zzz , \uuu  \Big)
&=\prod_{j=1}^{N} (1- z_j)^{- a } \oFz \Big(   \begin{matrix} a \\ \, \end{matrix} ;\frac{\zzz }{\one -\zzz } , \uuu-\one  \Big)
\\& \label{:23a}
=\prod_{j=1}^{N} (1- u_j)^{- a } \oFz \Big(   \begin{matrix} a \\ \, \end{matrix} ; \zzz-\one , \frac{\uuu }{\one -\uuu }  \Big)
,\\\label{:23b}
\zFz (; \zzz, \uuu)&=e^{ |\zzz |} \zFz (; \zzz, \uuu-\one)
.\end{align}
\begin{proof}
  
From the Jack binomial formula \eqref{:31f} and \eqref{:28a}, we obtain 
\begin{align*}\sum_{ \mm \in \PN }  (a  )_{\mm }^\theta \Psi_{\mm }^\theta (\zzz ) \Phi_{\mm}^\theta (\uuu )
& =
\sum_{ \mm, \kk \in \PN } \binom{\mm }{\kk }_{\theta } (a  )_{\mm }^\theta \Psi_{\mm }^\theta (\zzz ) \Phi_{\kk}^\theta (\uuu -\one )
\\&
=
 \prod_{j=1}^{N} (1-z_j)^{-a} \sum_{\kk \in \PN } (a)_{\kk}^\theta \Psi_{\kk }^\theta \Big(\frac{\zzz }{\one -\zzz } \Big) \Phi_{\kk}^\theta (\uuu  -\one )
,\end{align*}
which implies the first equality in \eqref{:23a}, whilst the second follows from the invariance of the left-hand side under the interchange of $\zzz$ and $\uuu$.
By replacing $\zzz $ by $ \zzz /a$ in \eqref{:23a} and letting $a \to \infty$, we obtain \eqref{:23b}.

\end{proof}
\end{lemma}

\section{Multivariate Meixner polynomials}\label{s:3}
We continue to assume that $\theta>0$ in this section.
We also retain the assumption that $(a)_{\mm}^{\theta}\neq0$ for all $\mm\in\PN$ and, in addition, assume that $c\notin \{ 0,1\}$ unless otherwise stated.

\subsection{Definition}
Define the multivariate Meixner polynomials by
\begin{align}\notag
M_{\mm }^{\theta, a , c} (\xxx  )
:=
M_{\mm }^{\theta} (\xxx  ; a , c)
:=\sum_{\kk \subseteq \mm } \frac{1}{(a )_{\kk } ^\theta \Psi_{\kk }^{\theta}( \one  ) } \binom{\mm }{\kk }_{\theta } \binom{\xxx }{\kk}_{ \theta } \Big(1-\frac{1}{c} \Big)^{|\kk |}   
.\end{align}
These Meixner polynomials were introduced and investigated for $\theta=1/2,1,2 $ in \cite{Shi16}.

Multivariate Meixner polynomials are often defined by a determinantal construction from the one-variable Meixner polynomials.
These Schur-type Meixner polynomials correspond to ours when $\theta=1$; see \sref{s:6}.
Note that other multivariate Meixner polynomials have also been studied in \cite{Gri75, Ili12}.

For $\mm, \xx\in \PN$, the duality 
\begin{align}\label{:31e}
M_{\mm}^{\theta,a,c}(\xx) =M_{\xx}^{\theta,a,c}(\mm ) 
\end{align}
follows immediately from the definition.
A direct computation yields the degeneration from the Meixner polynomials to the Laguerre polynomials: for any $\xxx \in \R_{}^N$, 
\begin{align}
 \lim_{c\to 1} M_{\mm }^{\theta,a , c}\Big(\frac{\xxx }{1-c}   \Big)
 =\frac{1}{(a)_{\mm}^\theta \Psi_{\mm} ^{\theta }(\one  )}  \LL_{\mm }^{\theta, a}(\xxx )
.\end{align}
Since $\{ P_{\mm }^{\theta, \ip }\}_{\mm\in \PN } $ forms a basis for $\Q(\theta )[\xxx ]^{\mathfrak{S}_N} $, the following claim is immediate.
\begin{proposition}\label{p:32}
The family $\{ M_{\mm }^{\theta, a , c} (\xxx -\theta \delta) \}_{\mm\in\PN }$ forms a basis for $\Q(\theta, a,c)[\xxx ]^{\mathfrak{S}_N} $.
\end{proposition}

\subsection{Generating functions }
We derive several types of double generating functions for multivariate Meixner polynomials.
We first show the following three double generating functions, which will later turn out to coincide with the double generating function of Laguerre polynomials.

\begin{theorem}\label{t:31}
For any $\gamma \in \R$, we have 
\begin{align} \notag&
 \sum_{\mm,\xx \in \PN }(\gamma )_{\mm}^\theta (a )_{\xx }^\theta  M_{\mm}^{\theta,a,c}(\xx) \Psi_{\mm}^\theta(c\uuu ) \Psi_{\xx }^\theta ( \zzz )
\\&\qquad   \label{:36b}
= \prod_{j=1}^{N}(1-c u_j)^{-\gamma }(1-z_j)^{-a } \oFz  \Big(\begin{matrix} \gamma \\  \, \end{matrix} \st  \frac{ \zzz }{\one -\zzz }, \frac{(c-1)\uuu }{\one -c \uuu } \Big) 
  \end{align}
  In particular, we have
  \begin{align}  \notag
  \sum_{\mm, \xx \in \PN  }  (a)_\mm^\theta  (a)_\xx^\theta M_{\mm}^{\theta,a,c}(\xx)  \Psi_{\mm }^\theta(c \uuu ) \Psi_{\xx }^\theta (\zzz ) 
=\prod_{j=1}^{N}(1- c u_j)^{-a} \oFz \Big(  \begin{matrix} a \\ \, \end{matrix} ;\zzz , \frac{\one -\uuu }{\one -c \uuu } \Big) 
\\  \label{:36a}
=\prod_{j=1}^{N}(1- z_j)^{-a} \oFz \Big(  \begin{matrix} a \\ \, \end{matrix} ; \frac{c \one -\zzz }{\one -\zzz},  \uuu \Big)
,  \\
\label{:36g}
 \sum_{\mm, \xx \in \PN } (a )_{\xx }^\theta  M_{\mm}^{\theta,a,c}(\xx) \Psi_{\mm}^\theta(c \uuu ) \Psi_{\xx }^\theta ( \zzz )
=
e^{c |\uuu |}\prod_{j=1}^{N}(1-z_j)^{-a } \zFz  \Big( \st  \frac{ \zzz }{\one -\zzz }, (c-1 )\uuu \Big) 
.\end{align}
\begin{proof}
We first show that the left-hand side of \eqref{:36b} converges absolutely for sufficiently small \(\|\uuu\|\) and $\| \zzz \|$.
The triangle inequality yields $  |(\gamma)_{\mm}^\theta | \le (|\gamma| +r)_{\mm}^\theta $, where $r=2\theta(N-1)$.
From this, we have
\begin{align*}&
 \sum_{\mm,\xx \in \PN } | (\gamma )_{\mm}^\theta  (a )_{\xx }^\theta  M_{\mm}^{\theta,a,c}(\xx) \Psi_{\mm}^\theta(c\uuu ) \Psi_{\xx }^\theta ( \zzz ) |
\\&
\le 
\! \sum_{\mm,\xx \in \PN } \!\! (|\gamma|+r )_{\mm}^\theta (|a| +r )_{\xx }^\theta \sum_{\kk \subseteq \mm } \frac{1}{|(a )_{\kk } ^\theta| \Psi_{\kk }^{\theta}( \one  ) } \binom{\mm }{\kk }_{\theta } \binom{\xx }{\kk}_{ \theta } \Big| 1-\frac{1}{c} \Big|^{|\kk |}  \Psi_{\mm}^\theta(|c| \| \uuu \| \one ) \Psi_{\xx }^\theta ( \| \zzz \| \one ) 
\\&
=
 \sum_{\kk \in \PN   } \frac{| 1-c^{-1} |^{|\kk |}}{|(a )_{\kk } ^\theta| \Psi_{\kk }^{\theta}( \one  ) } \frac{(|\gamma|+r )_{\kk }^\theta  (|a| +r)_{\kk }^\theta}{(1-|c| \|\uuu \| )^{(|\gamma |+r)N}(1-\|\zzz \| )^{(|a|+r)N}}  \Psi_{\kk }^\theta \Big(\frac{|c|\|\uuu\| }{1-|c|\|\uuu\| } \one \Big)  \Psi_{\kk }^\theta \Big(\frac{\|\zzz \| }{1-\|\zzz \| } \one \Big)
  \\&
= 
\frac{1}{(1-|c| \|\uuu \| )^{(|\gamma |+r )N}(1-\|\zzz \| )^{(|a|+r )N}} 
\\& \qquad\qquad\qquad
\times \sum_{\kk \in \PN   }  \frac{(|\gamma|+r )_{\kk }^\theta  (|a| +r )_{\kk }^\theta}{|(a )_{\kk } ^\theta| }   \Phi_{\kk }^\theta \Big(\frac{|c-1|\|\uuu\| }{1-|c|\|\uuu\| } \one \Big)  \Psi_{\kk }^\theta \Big(\frac{\|\zzz \| }{1-\|\zzz \| } \one \Big)
,\end{align*}
where we have used \eqref{:28a} in the first equality.
Thus, the absolute convergence of the left-hand side of \eqref{:36b} follows from the same argument as in \lref{l:B2}.

Equality \eqref{:36b} follows by using \eqref{:28a} twice.
Indeed, the left-hand side of \eqref{:36b} becomes
\begin{align*}&
\sum_{\mm ,\xx\in \PN  } \sum_{\kk \subseteq \mm } \frac{ (\gamma )_\mm ^\theta (a )_{\xx }^\theta}{(a )_{\kk }^\theta \Psi_{\kk}^\theta (\one )} \binom{\mm }{\kk }_{\theta } \binom{\xx  }{\kk}_{ \theta } \Big(1-\frac{1}{c} \Big)^{|\kk |}
\Psi_{\mm}^\theta(c \uuu ) \Psi_{\xx }^\theta (\zzz )
\\& 
= \sum_{\xx \in \PN  } \sum_{\kk \in \PN  }\frac{ (\gamma )_\kk^\theta (a )_{\xx }^\theta }{(a )_{\kk }^\theta }  \binom{\xx  }{\kk}_{ \theta } \Big(1-\frac{1}{c} \Big)^{|\kk |}
 \prod_{j=1}^{N} (1-cu_j)^{-\gamma } \Phi_{\kk }^\theta \Big(\frac{c \uuu }{\one  -c \uuu } \Big) \Psi_{\xx }^\theta ( \zzz )  \quad && \text{from \eqref{:28a}}
\\&
= \prod_{j=1}^{N} (1-c u_j)^{-\gamma }(1-z_j)^{-a }\sum_{\kk \in \PN  }  (\gamma )_\kk ^\theta   \Big(1-\frac{1}{c} \Big)^{|\kk |} \Phi_{\kk }^\theta \Big(\frac{c \uuu }{\one  -c \uuu } \Big) \Psi_{\kk }^\theta \Big(\frac{\zzz }{\one  -\zzz } \Big)   \quad && \text{from \eqref{:28a}}
.\end{align*}
Equality \eqref{:36a} follows from \eqref{:36b} with \eqref{:23a}.
Equality \eqref{:36g} follows from \eqref{:28a} and \eqref{:25b}, or degeneration from \eqref{:36b}.
\end{proof}
\end{theorem}





Hereafter, to simplify notation, we often use the symbol
\begin{align*}
  \ctil = \frac{c}{1-c} 
.\end{align*}
We have the following double generating functions, distinct from those in \tref{t:31}.
\begin{theorem} \label{t:32}
For any $\gamma \in \R$,  we have
\begin{align}\label{:27d}
\sum_{\mm, \xx \in \PN }(\gamma )_{\mm}^\theta  M_{\mm}^{\theta,a,c}(\xx)e^{-\ctil |\zzz | } \Psi_{\mm}^\theta(\uuu ) \Psi_{\xx }^\theta ( \ctil \zzz )
=
\prod_{j=1}^{N}(1-u_j)^{-\gamma } \oFo  \Big(\begin{matrix} \gamma \\ a  \end{matrix};-\zzz , \frac{\uuu }{\one -\uuu } \Big) 
.\end{align}
 In particular, we have
\begin{align} \label{:27c}
\sum_{\mm, \xx \in \PN  } (a)_\mm^\theta  M_{\mm}^{\theta,a,c}(\xx)e^{-\ctil |\zzz | } \Psi_{\mm }^\theta(\uuu ) \Psi_{\xx }^\theta (\ctil \zzz )
&=\prod_{j=1}^{N}(1-u_j)^{-a} \zFz  \Big(;-\zzz , \frac{\uuu }{\one -\uuu } \Big)
,\\
\label{:25a}
\sum_{\mm , \xx \in \PN  }M_{\mm}^{\theta,a,c}(\xx)e^{-\ctil |\zzz | } \Psi_{\mm}^\theta(\uuu ) \Psi_{\xx }^\theta ( \ctil \zzz )
&
=e^{|\uuu |} \zFo \Big(\begin{matrix}  \\ a  \end{matrix};-\zzz , \uuu  \Big) 
.\end{align}
\begin{proof}
Absolute convergence of the series on the left-hand side of \eqref{:27d} follows by an argument similar to that in the proof of \tref{t:31}, which also justifies the rearrangements below.
The left-hand side of \eqref{:27d} transforms into
\begin{align*}&
\sum_{\mm , \xx \in \PN  } \sum_{\kk \subseteq \mm } \frac{ (\gamma )_\mm ^\theta }{(a )_{\kk }^\theta \Psi_{\kk}^\theta (\one  )} \binom{\mm }{\kk }_{\theta } \binom{\xx  }{\kk}_{ \theta } \Big(1-\frac{1}{c} \Big)^{|\kk |}
e^{-\ctil |\zzz | } \Psi_{\mm}^\theta(\uuu ) \Psi_{\xx }^\theta (\ctil \zzz )
\\& 
= \sum_{\xx \in \PN  } \sum_{\kk \in \PN  }\frac{ (\gamma )_\kk^\theta }{(a )_{\kk }^\theta } \binom{\xx  }{\kk}_{ \theta } \Big(1-\frac{1}{c} \Big)^{|\kk |}  e^{-\ctil  |\zzz | } \prod_{j=1}^{N} (1-u_j)^{-\gamma } \Phi_{\kk }^\theta \Big(\frac{\uuu }{\one  -\uuu } \Big) \Psi_{\xx }^\theta (\ctil \zzz )  \quad && \text{from \eqref{:28a}}
\\&
= \prod_{j=1}^{N} (1-u_j)^{-\gamma }\sum_{\kk \in \PN  }  \frac{ (\gamma )_\kk ^\theta }{(a )_{\kk }^\theta }  \Big(1-\frac{1}{c} \Big)^{|\kk |} \Phi_{\kk }^\theta \Big(\frac{\uuu }{\one  -\uuu } \Big) \Psi_{\kk }^\theta (\ctil  \zzz )  \quad && \text{from \eqref{:25b}}
,\end{align*}
which is identical to the right-hand side of \eqref{:27d}.
Equation \eqref{:27c} is a specialisation of \eqref{:27d}; more precisely, it is equivalent to \eqref{:36g}, but we retain this formula for the sake of clarity.
Replacing $\uuu $ by $\uuu / \gamma $ and taking $\gamma \to\infty$ in \eqref{:27d}, we obtain \eqref{:25a}, or equation \eqref{:25a} can be proven by a direct computation with \eqref{:25b}.
Thus, we finish the proof.
\end{proof}
\end{theorem}

The right-hand side of \eqref{:27d}, \eqref{:27c}, and \eqref{:25a} coincides with those of \eqref{:26d}, \eqref{:26c}, and \eqref{:26a}, respectively.
Thus, the double generating functions for the Meixner polynomials are given by the generating functions for the Laguerre polynomials.
This observation plays a key role in the study of Meixner polynomials.
We derive generating functions from \tref{t:32}.

\begin{theorem}\label{t:36}
For all $\gamma \in \R$ and $\xx \in \PN $, we have
\begin{align}
 \label{:34b}
\sum_{\mm \in \PN }(\gamma)_{\mm}^\theta M_{\mm}^{\theta,a,c}(\xx ) \Psi_{\mm }^\theta(\uuu ) 
&= 
\prod_{j=1}^{N}(1-u_j)^{-\gamma } \sum_{\mm \in \PN }\frac{(\gamma )_\mm ^\theta}{(a)_\mm ^\theta}  \binom{\xx }{\mm}_\theta  \Phi_{\mm }^\theta \Big(\frac{c-1}{c} \frac{\uuu }{\one -\uuu } \Big)
.\end{align}
In particular,  we have
\begin{align}
\label{:34a}
\sum_{\mm \in \PN }  (a )_{\mm }^\theta M_{\mm } ^{\theta ,a,c} (\xx ) \Psi_{\mm }^{\theta } (\uuu )
&=
 \prod_{j=1}^{N}(1-u_j)^{-a }  \Phi_{\xx }^{\theta }\Big(\frac{\one -c^{-1} \uuu }{\one -\uuu } \Big)
,\\ \label{:27a}
\sum_{\mm \in \PN }M_{\mm }^{\theta,a,c}(\xx ) \Psi_{\mm }^\theta (\ctil \uuu )
&= \frac{e^{\ctil  |\uuu | } }{(a)_{\xx }^\theta \Psi_{\xx }^{\theta}( \one  ) }  \LL_{\xx }^{\theta,a}(\uuu ) 
.\end{align}
Here, identity \eqref{:27a} is valid for all $\uuu \in \R^N$.
\begin{proof}

From \eqref{:27d} and \eqref{:25b}, we have
\begin{align*}&
\sum_{\mm ,\xx\in \PN }(\gamma)_{\mm }^\theta  M_{\mm}^{\theta,a,c}(\xx ) \Psi_{\mm}^\theta(\uuu ) \Psi_{\xx }^\theta ( \ctil \zzz )
\\&=
e^{\ctil |\zzz | } \prod_{j=1}^{N}(1-u_j)^{-\gamma } \sum_{\mm \in \PN }\frac{(\gamma )_\mm ^\theta}{(a)_\mm ^\theta} \Psi_{\mm }^{\theta}(-\zzz  ) \Phi _{\mm }^\theta \Big( \frac{\uuu }{\one  -\uuu } \Big)
\\&
= \prod_{j=1}^{N}(1-u_j)^{-\gamma } \sum_{\mm \in \PN }\frac{(\gamma )_\mm ^\theta }{(a)_\mm ^\theta } (-\ctil )^{-|\mm |} \sum_{\xx \in \PN }\binom{\xx }{\mm }_\theta \Psi_{\xx }^\theta (\ctil \zzz )\Phi_{\mm }^\theta  \Big( \frac{\uuu }{\one  -\uuu } \Big)
.\end{align*}
Comparing the coefficients of $\Psi_{\xx }^\theta (\ctil \zzz )$ in this equation, we obtain \eqref{:34b}.
Equation \eqref{:34a} follows from \eqref{:34b} with \eqref{:31f}.
Replacing $\uuu $ by $\ctil \uuu /\gamma $ and taking $\gamma \to \infty $ in \eqref{:34b}, we obtain \eqref{:27a}.
\end{proof}
\end{theorem}

\begin{remark}
\begin{enumerate}
\item
For $\theta=1/2,1,2$, equation \eqref{:34a} was established in \cite[(3.22)]{Shi16} by another approach, which relies on symmetric cones and therefore does not apply to general $\theta$.
  

\item 
Another proof of \eqref{:27a}, which exhibits a direct connection between the Laguerre and Meixner polynomials, is obtained from the generating functions of the Laguerre polynomials \eqref{:26d}, \eqref{:26c}, or \eqref{:26a}.
Indeed, from \eqref{:25a} and \eqref{:26a}, we have 
\begin{align}\label{:34c}
\sum_{\mm \in \PN } \frac{1}{(a)_{\mm }^\theta } \LL_{\mm}^{\theta,a}(\zzz ) \Phi_{\mm }^\theta (\uuu )
=\sum_{\mm, \xx \in \PN  }M_{\mm}^{\theta, a,c}(\xx )e^{-\ctil |\zzz | } \Psi_{\mm}^\theta(\uuu ) \Psi_{\xx }^\theta (\ctil \zzz )
,\end{align}
which yields \eqref{:27a} by comparing the coefficients of $\Psi_{\mm}^\theta (\uuu )$.

\item 
When $N=1$, generating functions \eqref{:34b}, \eqref{:34a}, and \eqref{:27a} correspond to \cite[(9.10.13), (9.10.11), (9.10.12)]{KLS10}, respectively.
\end{enumerate}
\end{remark}

\subsection{Generating functions and generalised Laplace transforms}
We consider the generalised Laplace transform.
In this subsection, we assume that $a-\theta(N-1)>0$.
Define the normalising constant of the Laguerre ensemble, see \cite[Corollary 8.2.2]{AAR99}, by
\begin{align}\label{:34d}
\mathcal{N}_{\theta, a}= \int_{\R_{\ge}^N }  e^{-|\yyy | } \prod_{j=1}^N y_j^{a-1-\theta(N-1)} |\Delta(\yyy )|^{2\theta }  d\yyy 
 =\prod_{j=1}^{N}\frac{\Gamma(a-(N-j)\theta) \Gamma(1+j\theta )}{\Gamma(1+\theta)}
.\end{align}

\begin{lemma}\label{l:37}
Assume that $p\le q$.
Then, for $\zzz \in (0,\infty)^N $, we have
\begin{align}\notag
  \frac{1}{\mathcal{N}_{\theta ,a}} \int_{\R_{\ge}^N }  \pFq \Big(\begin{matrix} a_1 , \ldots, a_p  \\ b_1 , \ldots , b_q  \end{matrix} \st \yyy, \uuu \Big) \zFz (; -\yyy,\zzz ) \prod_{k=1}^N y_k^{a-1-\theta(N-1)} |\Delta (\yyy )|^{2\theta } d\yyy
\\\label{:36e}
  =\Big( \prod_{j=1}^{N}z_j^{-a} \Big)\,\, {{}^{}_{p+1} F _{q} ^\theta } \Big(\begin{matrix} a_1 , \ldots, a_p ,a \\ b_1 , \ldots , b_q  \end{matrix} \st \frac{1}{\zzz }, \uuu \Big) 
.\end{align}
\begin{proof}
We first show that the integral in \eqref{:36e} is well-defined.
Let $M= \max_{i} z_i $ and $m= \min_{i} z_i $.
From \lref{l:1}, \eqref{B:1a}, and \eqref{:25b}, we have $$ \zFz (; \yyy, M \one -\zzz  ) \le \sum_{\mm \in \PN} \Psi_{\mm}^\theta (\yyy) (M-m)^{|\mm |} =e^{(M-m) |\yyy |}.$$
Hence, by \eqref{:23b}, we obtain $
  \zFz (; -\yyy, \zzz )
  = e^{-M | \yyy| } \zFz (; \yyy, M \one -\zzz  ) 
  \le e ^{- m |\yyy | }
$, which, together with \eqref{B:2c}, proves that the integral on the left-hand side of \eqref{:36e} converges absolutely for sufficiently small $\|\uuu \|$.
Using the formula
\begin{align}\label{:36f}
\frac{1}{\mathcal{N}_{\theta ,a}}  \int_{\R_{\ge}^N } \Phi_{\kk}^\theta (\yyy ) \, \zFz (; -\yyy,\zzz ) \prod_{j=1}^N y_j^{a-1-\theta(N-1)} |\Delta (\yyy )|^{2\theta } d\yyy
= (a)_{\kk}^\theta \Big( \prod_{j=1}^{N}z_j^{-a}  \Big)
\Phi_{\kk}^\theta \Big( \frac{1}{\zzz } \Big)
\end{align}
given by \cite[(6.1)]{BaF97} (note that it has a small typo), a direct computation yields \eqref{:36e}. 
\end{proof}
\end{lemma}

\begin{theorem}\label{t:37}
Assume that $a-\theta(N-1)>0$.
Then, for all $\gamma \in \R $, we have 
\begin{align}& \notag
\frac{1}{\mathcal{N}_{\theta ,a}} \int_{\R_{\ge}^N }  \prod_{j=1}^{N}(1-u_j)^{-\gamma } \oFo \Big( \begin{matrix}
 \gamma \\ a  
 \end{matrix} ;
 (c-1 ) \yyy, \frac{\uuu }{\one -\uuu } \Big) \zFz (;-\yyy , \one -\zzz) \prod_{j=1}^N y_j^{a-1-\theta(N-1)} |\Delta (\yyy )|^{2\theta } d\yyy 
\\ & \label{:39f}
=\prod_{j=1}^{N}(1-c u_j)^{-\gamma}(1- z_j)^{- a } \oFz \Big(
\begin{matrix} \gamma \\ \, \end{matrix} ;\frac{\zzz }{\one -\zzz }, \frac{(c-1)\uuu }{\one -c\uuu } \Big)
.\end{align}
In particular, we have
\begin{align}
& \notag
 \frac{1}{\mathcal{N}_{\theta, a}} \int_{\R_{\ge}^N }    \prod_{j=1}^{N}(1-u_j)^{-a } \zFz \Big( ;(c -1) \yyy , \frac{\uuu }{\one -\uuu } \Big)  \zFz (;-\yyy , \one -\zzz ) 
 \prod_{k=1}^N y_k^{a-1-\theta(N-1)} |\Delta (\yyy )|^{2\theta } d\yyy 
\\ & \label{:39g}
=\prod_{j=1}^{N}(1-c u_j)^{-a}  \oFz \Big(
\begin{matrix} a \\ \, \end{matrix} ; \zzz , \frac{\one  - \uuu }{\one -c\uuu } \Big)
,\\& \notag
 \frac{1}{\mathcal{N}_{\theta, a}} \int_{\R_{\ge}^N } e^{|\uuu |} \zFo \Big(\begin{matrix}
  \\
  a
 \end{matrix} ;(c -1) \yyy , \uuu \Big)
 \zFz (; -\yyy , \one -\zzz )  \prod_{k=1}^N y_k^{a-1-\theta(N-1)} |\Delta (\yyy )|^{2\theta } d\yyy 
\\ & \label{:39h}
=e^{c|\uuu | }\prod_{j=1}^{N}(1-z_j)^{-a}  \zFz \Big(; \frac{\zzz }{\one -\zzz } , (c-1)\uuu  \Big)
.\end{align}
  \begin{proof}
Direct computations using \eqref{:36e} and \eqref{:23a} derive \eqref{:39f}. 
Equations \eqref{:39g} and \eqref{:39h} follow from \eqref{:39f}.
\end{proof}
\end{theorem}

\tref{t:37} shows that the double generating functions for the Laguerre polynomials in \eqref{:39f}, \eqref{:39g}, and \eqref{:39h} correspond, respectively, to the double generating functions for the Meixner polynomials in \eqref{:36b}, \eqref{:36a}, and \eqref{:36g}.
In view of this correspondence, the following theorem shows that a generalised Laplace transform 
corresponds to a generating function for the Meixner polynomials.

\begin{theorem}\label{t:35}
Assume that  $a-\theta(N-1)>0$ holds.
Then, for all $\mm \in \PN $, we have
\begin{align}\notag&
\frac{1}{\mathcal{N}_{\theta, a}} \int_{\R_{\ge}^N }   \LL_{ \mm}^{\theta, a} ((1-c) \yyy ) \zFz (; -\yyy , -\zzz )e^{-|\yyy |} \prod_{j=1}^N y_j^{a-1-\theta(N-1)} |\Delta(\yyy )|^{2\theta } d\yyy  
\\&  \label{:36c} \qquad 
=(a)_{\mm}^\theta\Psi_{\mm }^{\theta}( c \one ) \sum_{\xx \in \PN }  (a )_{\xx }^\theta  M_{\xx } ^{\theta ,a,c } (\mm ) \Psi_{\xx }^{\theta } (\zzz )
.\end{align}
In particular, for every $\xx\in \PN $, we have 
\begin{align}\notag 
\frac{1}{\mathcal{N}_{\theta, a} }  \int_{\R_{\ge}^N } \LL_{\mm}^{\theta,a} ((1-c) \yyy ) \frac{\Phi_\xx^{\theta} (\yyy ) }{(a)_{\xx }^\theta } e^{-|\yyy | } \prod_{j=1}^N y_j^{a-1-\theta(N-1)} |\Delta(\yyy )|^{2\theta }  d\yyy
\\ \label{:36d} \qquad   
= (a )_{\mm }^\theta \Psi_{\mm}^{\theta}(c \one  ) M_{\mm } ^{\theta ,a,c} (\xx ) 
.\end{align}
Furthermore, letting $c\to 0$ in \eqref{:36d}, we have 
\begin{align*}
\frac{1}{\mathcal{N}_{\theta, a} } \int_{\R_{\ge}^N } \LL_{\mm}^{\theta,a} ( \yyy ) \frac{\Phi_\xx ^{\theta} (\yyy ) }{(a)_{\xx }^\theta } e^{-|\yyy | } \prod_{j=1}^N y_j^{a-1-\theta(N-1)} |\Delta(\yyy )|^{2\theta } d\yyy 
= (-1)^{|\mm |} \binom{\xx } {\mm }_\theta
.\end{align*}
\begin{proof}

From \eqref{:26d}, \eqref{:36b}, and \eqref{:39f}, it follows that 
\begin{align*}&
  \sum_{\mm,\xx \in \PN }(\gamma )_{\mm}^\theta (a )_{\xx }^\theta  M_{\mm}^{\theta,a,c}(\xx) \Psi_{\mm}^\theta(c\uuu ) \Psi_{\xx }^\theta ( \zzz )
   \\&
  =\frac{1}{\mathcal{N}_{\theta ,a}} \int_{\R_{\ge}^N } \sum_{\mm \in \PN } \frac{(\gamma )_{\mm} ^\theta }{(a)_{\mm} ^\theta }\LL_{\mm}^{\theta,a}((1-c)\yyy ) \Phi_{\mm }^\theta (\uuu ) \zFz (;-\yyy , \one -\zzz) \prod_{k=1}^N y_k^{a-1-\theta(N-1)} |\Delta (\yyy )|^{2\theta } d\yyy 
.\end{align*}
Thus, comparing the coefficients of $\Phi_\mm^\theta (\uuu )$ in this equation, we have \eqref{:36c}.
Comparing the coefficients of $\Phi_{\xx}^\theta (\zzz )$ in \eqref{:36c} and by the duality \eqref{:31e}, we obtain \eqref{:36d}.
\end{proof}  
\end{theorem}

\begin{remark}
\begin{enumerate}
  \item 
Equation \eqref{:36c} follows from \eqref{:26c}, \eqref{:36a}, and \eqref{:39g}, or \eqref{:26a}, \eqref{:36g}, and \eqref{:39h} instead of \eqref{:26d}, \eqref{:36b}, and \eqref{:39f}.

  \item
Equation \eqref{:36c} can be shown by a direct computation based on the equation 
\begin{align}\notag
\frac{1}{\mathcal{N}_{\theta, a}} \int_{\R_{\ge}^N } \LL_{ \mm }^{\theta,a } (\yyy ) \zFz (;-\yyy , \uuu ) \prod_{j=1}^N y_j^{a-1-\theta(N-1)} | \Delta(\yyy )|^{2\theta } d\yyy   
\\ \label{:35a}
=(a)_{\mm}^\theta \Big( \prod_{j=1}^N u_j^{-a } \Big) \Psi_{\mm}^{\theta} \Big( \frac{\uuu -\one }{\uuu } \Big)
\end{align}
established in \cite[(6.2)]{BaF97}.
Actually, by the change of variables $\yyy \mapsto (1-c)\yyy $ and $\uuu \mapsto (\one -\uuu )/(1-c)$ in \eqref{:35a}, we have
\begin{align*}&
\frac{1}{\mathcal{N}_{\theta, a}} \int_{\R_{\ge}^N } \LL_{ \mm}^{\theta, a} ((1-c) \yyy ) \zFz (; -\yyy, -\uuu ) e^{-|\yyy |}\prod_{j=1}^N y_j^{a-1-\theta(N-1)} |\Delta(\yyy )|^{2\theta } d\yyy 
\\&
= (a)_{\mm}^\theta c^{|\mm|} \prod_{j=1}^N (1-u_j)^{-a } \Psi_{\mm}^{\theta} \Big( \frac{\one -c^{-1}\uuu }{\one -\uuu } \Big)
\\&
=(a)_{\mm}^\theta  \Psi_{\mm}^{\theta}( c \one ) \sum_{\kk \in \PN }  (a )_{\kk }^\theta M_{\kk } ^{\theta ,a , c} (\mm ) \Psi_{\kk }^{\theta } (\uuu )
.\end{align*}
Here, we have used \eqref{:34a} and \eqref{:23b}.
\end{enumerate}  
\end{remark}

\subsection{Difference equations }
This subsection is devoted to deriving a difference operator having the Meixner polynomials as eigenfunctions.
To this end, we first derive their difference equations.

We define $h_{+ ,i}^\theta$ and $h_{- ,i}^\theta$ by 
\begin{align*}
h_{\pm ,i}^\theta (\xxx )=\prod_{1\le j \neq i \le N } \frac{x_i-x_j-\theta(i-j)\pm \theta }{x_i-x_j-\theta(i-j) } 
,\end{align*}
where the same sign is chosen on both sides.
We then have 
\begin{align}\label{:315c}
\binom{\mm}{ \mm-\eee _{i}} _\theta 
= (m_i +\theta (N-i)) h_{-,i}^\theta(\mm ) 
= \frac{\Psi _{\mm-\eee _i }^\theta (\one  )}{\Psi _{\mm }^\theta (\one )} h_{+,i}^\theta(\mm -\eee _i)
.\end{align}
Here, the first equality follows from \cite[Section 2.5.4, (1)]{ViK95}, whilst the second is a consequence of the following identity, obtained directly from \eqref{:20a} and \eqref{:20b}:
\begin{align}\label{:315a}
\frac{h_{+,i}^\theta(\mm ) }{h_{-,i}^\theta(\mm +\eee _i )} =(m_i+1+\theta(N-i) ) \frac{\Psi_{\mm +\eee_ i}^\theta( \one )}{\Psi_{\mm }^\theta (\one  )} 
.\end{align}

We quote the following actions from \cite{Kan93}:
\begin{align}\label{:315g}
E_0^{N} \Psi_{\mm}^\theta (\xxx )&= \Psi_{\mm}^\theta (\one  ) \sum_{i=1}^{ N } \binom{\mm}{ \mm-\eee _{i}} _\theta \Phi_{\mm-\eee _{i}}^\theta  (\xxx )
,\\ \label{:315h}
E_1^{N} \Psi_{\mm}^\theta (\xxx )&=|\mm | \Psi_{\mm}^\theta (\xxx ) 
,\\ \label{:315j}
D_1^{N,\theta} \Psi_{\mm}^\theta (\xxx )&= \Psi_{\mm}^\theta (\one  )  \sum_{i=1}^{ N } \binom{\mm }{ \mm -\eee _{i}} _\theta  (m_i-1+\theta(N-i) )  \Phi_{\mm-\eee _{i}}^\theta  (\xxx  )
,\\ \label{:315k}
|\xxx | \Psi_{\mm }^\theta (\xxx )&=\sum_{i=1}^{ N } \binom{\mm +\eee _i }{ \mm }_\theta \Psi_{\mm +\eee _i }^\theta (\xxx )
.\end{align}
The summations in \eqref{:315g} and \eqref{:315j} are over all $i$ such that $\mm -\eee_i $ is admissible, and the summation in \eqref{:315k} is over all $i$ such that $\mm +\eee_i $ is admissible.
This convention will be used in all future summations involving $\mm \pm \eee_i $ unless otherwise stated.

\begin{lemma} \label{l:31}
For all $\mm \in \PN $ and $\xxx \in\R^N $, we have
\begin{align}\label{:315d}
E_0^{N} \Psi_{\mm}^\theta (\ctil \xxx )&= \ctil \sum_{i=1}^{ N }h_{+,i}^{\theta } (\mm -\eee _i) \Psi_{\mm - \eee _i}^\theta  (\ctil \xxx )
 ,\\\label{:315e}
D_1^{N,\theta } \Psi_{\mm}^\theta (\ctil \xxx )&= \ctil \sum_{i=1}^{ N }   (m_i-1+\theta (N-i)) h_{+,i}^{\theta} (\mm -\eee _i) \Psi_{\mm -\eee _i}^\theta (\ctil \xxx )
,\\ \label{:315f}
\ctil |\xxx | \Psi_{\mm}^\theta (\ctil \xxx )&= \sum_{i=1}^{ N } (m_i+1+\theta (N-i)) h_{-,i}^{\theta} (\mm +\eee _i) \Psi_{\mm +\eee_i}^\theta (\ctil \xxx )
,\end{align}
\begin{proof}

Combining \eqref{:315g}, \eqref{:315j}, and \eqref{:315k} with \eqref{:315c}, we obtain the desired results.
\end{proof}  
\end{lemma}

The difference equation follows from the same computation as in \cite{Shi16}, but we include the proof for the reader’s convenience.
Note that although \pref{p:39} and \pref{p:312} hold for $\xxx \in \mathbb{R}^N$, it suffices to verify them for $\xx \in \PN$, since both sides are rational functions.

\begin{proposition}\label{p:39}
For all $\mm \in \PN$ and $\xxx \in \R^N$, we have
  \begin{align}\notag 
  (c-1)|\mm| M_{\mm}^{\theta,a,c}(\xxx )
&
  =c \sum_{i=1}^N (x_i+a-\theta(i-1))h_{+,i}^\theta (\xxx ) M_{\mm}^{\theta, a,c}(\xxx +\eee_i )
  \\ &  \notag \quad
  -((1+c)|\xxx | +acN )M_{\mm}^{\theta,a,c}(\xxx ) 
  \\&   \label{:311a}\quad
  +\sum_{i=1}^N (x_i+\theta(N-i))h_{-,i}^\theta (\xxx ) M_{\mm}^{\theta,a,c}(\xxx -\eee _i )
.\end{align}
\begin{proof}
Applying $\mathsf{A}_{\theta, a}^{N, (\zzz )} $ to \eqref{:26c} and \eqref{:27c}, we have 
\begin{align}\label{:31b}
  \sum_{\mm \in \PN } \mathsf{A}_{\theta, a }^{N, (\zzz )}  \LL_{\mm}^{\theta,a}(\zzz ) \Phi_{\mm }^\theta (\uuu )
  =  \sum_{\mm, \xx \in \PN  } (a)_\mm^\theta M_{\mm}^{\theta,a,c}(\xx ) \Psi_{\mm}^\theta(\uuu ) \mathsf{A}_{\theta, a }^{N, (\zzz )} ( e^{-\ctil |\zzz | } \Psi_{\xx }^\theta (\ctil \zzz ) )
.\end{align}
From \eqref{:31a} and \eqref{:315h}, the left-hand side of \eqref{:31b} transforms into 
\begin{align} \notag
-\sum_{\mm \in \PN } |\mm| \LL_{\mm}^{\theta,a}(\zzz ) \Phi_{\mm }^\theta (\uuu ) 
 & 
 =
-\sum_{\mm \in \PN } \LL_{\mm}^{\theta,a}(\zzz ) E_{1 }^{N,(\uuu )}\Phi_{\mm }^\theta (\uuu ) 
 \\&\notag
=
-E_{1}^{N,(\uuu )} \sum_{\mm,\xx \in \PN  } (a)_\mm^\theta M_{\mm}^{\theta,a,c}(\xx )e^{-\ctil |\zzz | } \Psi_{\mm}^\theta(\uuu ) \Psi_{\xx }^\theta (\ctil \zzz )
\\& \label{:31d}
=
-\sum_{\mm,\xx \in \PN  } |\mm| (a)_\mm^\theta M_{\mm}^{\theta,a,c}(\xx )e^{-\ctil |\zzz | } \Psi_{\mm}^\theta(\uuu ) \Psi_{\xx }^\theta (\ctil \zzz )
.\end{align}

Since $\mathsf{A}_{\theta, a}^{N, (\zzz )}$ is a second differential operator, the Baker-Campbell-Hausdorff formula yields
\begin{align}\label{:39a}
e^{\ctil |\zzz |} \mathsf{A}_{\theta, a}^{N, (\zzz )} e^{-\ctil |\zzz |}
 = e^{ \mathrm{ad } (\ctil |\zzz | )} \mathsf{A}_{\theta, a}^{N, (\zzz )}
 =\sum_{k=0}^{2} \frac{\mathrm{ad } (\ctil |\zzz | ) ^k}{k!} \mathsf{A}_{\theta, a}^{N, (\zzz )}
.\end{align}
Furthermore, the following relations hold from \cite{DeH12}:
\begin{align*}
[|\zzz |, D_{1}^{N, (\zzz )} ]=-2E_{1}^{N,(\zzz )} -\theta N(N-1)
, \quad [| \zzz |, E_{1}^{N,(\zzz )}  ]=-|\zzz |
, \quad   [|\zzz |, E_{0}^{N, (\zzz )} ]=-N
.\end{align*}
Therefore, we have 
\begin{align}\label{:39b}
\mathrm{ad } (\ctil |\zzz | ) \mathsf{A}_{\theta, a}^{N,(\zzz )}
=\ctil (-2E_{1}^{N,(\zzz )} + |\zzz | -aN )
,
\qquad 
\frac{\mathrm{ad } (\ctil |\zzz | ) ^2}{2!} \mathsf{A}_{\theta, a}^{N,(\zzz )}
=\ctil ^2 |\zzz |
.\end{align}
Substituting \eqref{:39b} into \eqref{:39a}, and using \eqref{:315h},  \eqref{:315d}, \eqref{:315e}, and \eqref{:315f}, we obtain 
\begin{align}& \notag
 e^{ \ad (\ctil |\zzz | )} \mathsf{A}_{\theta, a}^{N,(\zzz )} \Psi_{\xx }^\theta (\ctil \zzz )
 \\& \notag
 =\big\{ D_{1}^{N ,(\zzz )} - \frac{1+c}{1-c}E_{1}^{N, (\zzz ) }+(a-\theta(N-1)) E_0^{N,(\zzz )}+ \frac{\ctil }{1-c} |\zzz | - a \ctil N \big\}\Psi_{\xx }^\theta (\ctil \zzz )
 \\& \notag 
=  \ctil \sum_{i=1}^{N } (x_i-1+a-\theta (i-1)) h_{+,i}^{\theta} (\xx -\eee _i) \Psi_{\xx -\eee _i}^\theta (\ctil \zzz  )- \Big( \frac{1+c}{1-c} |\xx|+ a \ctil N  \Big)\Psi_{\xx }^\theta (\ctil \zzz  )
 \\&  \label{:39d}\quad
 +\frac{1 }{1-c}\sum_{i=1}^{N } (x_i+1+\theta (N-i)) h_{-,i}^{\theta} (\xx +\eee _i) \Psi_{\xx +\eee _i}^\theta (\ctil \zzz ) 
.\end{align}
Hence, the right-hand side of \eqref{:31b} transforms into
\begin{align}& \notag
  \sum_{\mm , \xx \in \PN  } (a)_\mm^\theta  M_{\mm}^{\theta,a,c}(\xx )  e^{-\ctil |\zzz | }  \Psi_{\mm}^\theta(\uuu ) 
 \\& \quad \notag 
\times \Big\{ \ctil \sum_{i=1}^{N } (x_i-1+a-\theta (i-1)) h_{+,i}^{\theta} (\xx -\eee _i) \Psi_{\xx -\eee _i}^\theta (\ctil \zzz  )- \Big( \frac{1+c}{1-c} |\xx |+a \ctil  N  \Big)\Psi_{\xx }^\theta (\ctil \zzz )
 \\& \quad \qquad \label{:31c}
 +\frac{1 }{1-c}\sum_{i=1}^{N } (x_i+1+\theta (N-i)) h_{-,i}^{\theta} (\xx +\eee _i) \Psi_{\xx +\eee _i}^\theta (\ctil \zzz ) \Big\}
.\end{align}
Comparing the coefficients of $\Psi_\mm^{\theta}(\uuu )\Psi_{\xx }^\theta (\ctil \zzz )$ on both sides of \eqref{:31b}, using \eqref{:31d} and \eqref{:31c}, we obtain \eqref{:311a} for $\xx \in \PN $.
Here, the summation is taken over all $i$, since the corresponding coefficient vanishes whenever the associated index is not admissible.
\end{proof}
\end{proposition}

\begin{remark}
By applying the same argument as in \pref{p:39} to \eqref{:27a}, we obtain  the recurrence relation in the index for the multivariate Meixner polynomials. 
However, by the duality \eqref{:31e}, the recurrence relation is equivalent to the difference equation \eqref{:311a}.
\end{remark}

We next rewrite \eqref{:311a} in terms of difference operators.
\begin{proposition}\label{p:312}
For all $\xxx \in \R^N $, we have
\begin{align}\label{:33e}
\sum_{i=1}^N (x_i-\theta(i-1))  \prod_{1\le j \neq i \le N } \frac{x_i-x_j-\theta(i-j-1) }{x_i-x_j-\theta(i-j) }   &= |\xxx | 
,\\ \label{:33f}
\sum_{i=1}^N  (x_i+\theta(N-i))\prod_{1\le j \neq i \le N } \frac{x_i-x_j-\theta(i-j+1)}{x_i-x_j-\theta(i-j) }  &=|\xxx | 
.\end{align}
In particular, we have
\begin{align}\notag &
 \sum_{i=1}^N \Big \{ \frac{c(x_i+a-\theta(i-1))}{1-c}  \prod_{1\le j \neq i \le N } \frac{x_i-x_j-\theta(i-j-1) }{x_i-x_j-\theta(i-j) }  \nabla _{+,i }  
 \\ \notag & \qquad \qquad
 +\frac{x_i+\theta(N-i) }{1-c}\prod_{1\le j \neq i \le N } \frac{x_i-x_j-\theta(i-j+1)}{x_i-x_j-\theta(i-j) }  \nabla_{-,i}\Big\} M_{\mm}^{\theta,a,c}(\xxx )
\\  \label{:33c} &
 = -|\mm| M_{\mm}^{\theta,a,c}(\xxx )
.\end{align}
\begin{proof}
We write $s_i:=x_i- \theta i $ and $\sss =(s_1,\ldots, s_N)$ to simplify notation. 
%
By applying the Lagrange interpolation formula, we have 
\begin{align*}
\prod_{j=1}^{N} (u- s_j \pm \theta )- \prod_{j=1}^{N} (u-s_j) 
&
=\sum_{i=1}^N \prod_{j=1}^{N} (s_i- s_j \pm \theta )\prod_{j \neq i }^{N} \frac{u-s_j}{s_i-s_j}
\\&
=\pm \theta \sum_{i=1}^N  \prod_{j \neq i }^{N} \frac{s_i- s_j \pm \theta }{s_i-s_j}(u-s_j)
,\end{align*}
where the signs are taken consistently.
Hence, it holds that 
\begin{align}\label{:33a}
 \prod_{j=1}^{N} \frac{u- s_j \pm \theta }{ u-s_j}=1\pm \theta \sum_{i=1}^N \frac{1 }{u-s_i}\prod_{j \neq i }^{N} \frac{s_i- s_j \pm \theta }{s_i-s_j}
.\end{align}

Substituting $u=v^{-1}$ in \eqref{:33a}, we have
\begin{align}\notag 
 \prod_{j=1}^{N} \frac{1- s_j v \pm \theta v}{ 1-s_j v}
 &=1\pm \theta \sum_{i=1}^N \frac{ v }{1-s_i v}\prod_{j \neq i }^{N} \frac{s_i- s_j \pm \theta }{s_i-s_j}
 \\& \label{:313b}
 =1 \pm \theta \sum_{i=1}^N \sum_{k\ge 0 } s_i^k v^{k+1}\prod_{j \neq i }^{N} \frac{s_i- s_j \pm \theta }{s_i-s_j}
.\end{align}
Expanding the left-hand side of \eqref{:313b} in another way, we have 
\begin{align}\notag 
  \prod_{j=1}^{N} \frac{1- s_j v \pm \theta v}{ 1-s_j v}
  &
  =\sum_{k\ge 0} h_{k}(\sss ) v^k  \sum_{l\ge 0} e_{l} (\pm \theta \one  -\sss )v^l
  \\& \label{:313c}
  =\sum_{m\ge 0}  \sum_{l= 0}^{m}  h_{m-l}(\sss ) e_{l} (\pm \theta \one  -\sss )v^m
,\end{align}
where $e_k$ and $h_k$ are the elementary, and complete homogeneous symmetric polynomials, respectively.
Comparing the coefficients of $v$, $ v^2$ in the expansions \eqref{:313b} and \eqref{:313c}, we have  
\begin{align*}
\pm \theta \sum_{i=1}^N \prod_{j \neq i }^{N} \frac{s_i- s_j \pm \theta }{s_i-s_j}
=\pm \theta N
,\qquad 
\pm \theta \sum_{i=1}^N s_i \prod_{j \neq i }^{N} \frac{s_i- s_j \pm \theta }{s_i-s_j}
=\pm \theta |\sss |+\frac{N(N-1)}{2} \theta ^2
,\end{align*}
which imply \eqref{:33e} and \eqref{:33f}.
Equation \eqref{:33c} immediately follows from \pref{p:39} together with \eqref{:33e} and \eqref{:33f}.
Thus, the proof is complete.
\end{proof}
\end{proposition}

\subsection{Orthogonality}
The Meixner polynomials in one variable are orthogonal with respect to the negative binomial distribution.
This orthogonality extends to the multivariate case.

\begin{proposition}\label{p:313}
Assume that $0<c<1$.
Then, for any $\mm, \kk\in\PN $, we have 
\begin{align}\label{:33d}
  \sum_{\xx \in \PN } (a)_{\xx}^\theta (1-c)^{Na}  \Psi_{\xx }^\theta (c \one ) M_{\mm}^{\theta,a,c}(\xx )M_{\kk}^{\theta,a,c}(\xx) = \frac{1}{ (a)_{\mm}^\theta \Psi_{\mm}^{\theta}(c \one )}\delta_{\mm,\kk }
.\end{align}
Here, the orthogonality measure is positive if $a-\theta(N-1)>0$. 
\begin{proof}
 Taking $\zzz =(c \one -\www ) / (\one -\www )$ in \eqref{:36a}, and using  \eqref{:34a}, we have 
  \begin{align*}&
\sum_{\mm,\kk \in \PN } (a)_{\mm}^\theta \Phi_{\mm}^\theta(\www )\Psi_{\kk}^\theta (\uuu )\delta_{\mm,\kk}
    \\&
    =\oFz   \Big(
  \begin{matrix} a \\ \, \end{matrix} ;\www , \uuu \Big)
\\&
= (1-c)^{Na} \prod_{j=1}^{N}(1-w_j)^{-a} \sum_{\mm ,\xx\in \PN  } (a)_\mm^\theta (a)_\xx^\theta M_{\mm}^{\theta, a, c}(\xx ) \Psi_{\mm }^\theta (c\uuu ) \Psi_{\xx }^\theta \Big(\frac{c \one -\www }{\one -\www } \Big) 
\\&
= (1-c)^{Na}  \sum_{\mm, \xx \in \PN  } (a)_\mm^\theta (a)_\xx ^\theta M_{\mm}^{\theta,a, c}(\xx) \Psi_{\mm }^\theta (c \uuu ) \Psi_{\xx}^\theta (c \one  )  \sum_{\kk \in \PN }  (a )_{\kk }^\theta M_{\kk } ^{\theta, a, c } (\xx ) \Psi_{\kk }^{\theta } ( \www )
.\end{align*}
Comparing the coefficients of $\Phi_{\mm}^\theta(\uuu )\Psi_{\kk}^\theta (\www )$ on both sides, we obtain the desired result.
The estimate in the proof of \tref{t:31} justifies the rearrangements and coefficient comparisons for sufficiently small $\| \uuu \| ,\|\www\|$.
\end{proof}
\end{proposition}

\begin{remark}
\pref{p:313} is derived also from the kernel function for Laguerre  
\begin{align*}
\sum_{\xx\in \PN } \frac{c^{|\xx |} }{(a)_{\xx }^\theta \Psi_{\xx }^{\theta}( \one ) }  \LL_{\xx }^{\theta ,a } \Big( \frac{\zzz }{\ctil} \Big) \LL_{\xx }^{\theta ,a}\Big( \frac{\www }{\ctil} \Big) =(1-c)^{-Na} e^{-|\zzz | } e^{-|\www | } \zFo \Big( \begin{matrix} \\ a \end{matrix} \st  \frac{\zzz }{c} , \www \Big)
\end{align*}
introduced in \cite[Proposition 4.12]{BaF97}, and vice versa.
Indeed, using \eqref{:27a}, we have
\begin{align*}
(1-c)^{-Na}  \zFo \Big(\begin{matrix} \\ a \end{matrix}; \frac{\zzz }{c} , \uuu \Big)
&=
\sum_{\xx\in \PN } \frac{c^{|\xx |} }{(a)_{\xx }^\theta \Psi_{\xx }^{\theta}( \one ) }  e^{|\zzz | } \LL_{\xx }^{\theta ,a } \Big(\frac{ \zzz }{\ctil } \Big) e^{|\uuu | }  \LL _{\xx }^{\theta ,a}\Big(\frac{ \uuu }{\ctil } \Big)
\\&
=\sum_{\mm, \kk  \in \PN } \sum_{\xx\in \PN } c^{|\xx |} (a)_{\xx }^\theta \Psi_{\xx }^{\theta}( \one ) M_{\mm }^{\theta,a, c}(\xx )  M_{\kk }^{\theta, a, c}(\xx ) \Psi_{\mm }^\theta (\zzz ) \Psi_{\kk }^\theta ( \uuu ) 
.\end{align*}
Therefore, the orthogonality follows by comparing the coefficients of $\Psi_{\mm }^\theta (\zzz ) \Psi_{\kk }^\theta ( \uuu ) $.
\end{remark}

\pref{p:313} with \eqref{:27a} and \eqref{:31e} yields the following kernel function.
\begin{proposition}
Assume that $0<c<1$.
Then, for any $\xx \in \PN $ and $\uuu \in \R^N $, we have
\begin{align}\label{:318a}
  \sum_{\mm \in \PN }^{} M_{\mm}^{\theta, a,c}(\xx ) \LL_{\mm}^{\theta,a }(\uuu ) c ^{|\mm |} 
  =\frac{(1-c)^{-|\xx |-Na }e^{-\ctil |\uuu |}}{(a)_\xx^\theta}  \Phi_{\xx }^\theta (\uuu)
.\end{align}
\end{proposition}

\begin{remark}
When $a-\theta(N-1) >0$, the orthogonality of the Laguerre and Meixner polynomials shows that the three identities \eqref{:27a}, \eqref{:36d}, and \eqref{:318a} are equivalent: each of them implies the other two.
To the best of our knowledge, the identity \eqref{:318a} is not known even for the case $N=1$.
\end{remark}

\subsection{Inversion formulae}
Recall that the Meixner polynomials are defined as linear combinations of the interpolation Jack polynomials.
Here, the converse statement is established: the interpolation Jack polynomials can be expressed as linear combinations of the Meixner polynomials; see \cite[Section~6]{KoS98} for the one-variable case.
We first prepare an auxiliary lemma.

\begin{lemma}\label{l:314}
For a family of real numbers $\{a_\mm\}_{\mm \in \PN }$, let $\{ b_{\mm}\}_{\mm \in \PN } $ be the binomial transform of $\{a_\mm\}_{\mm \in \PN }$ defined by 
\begin{align*}
  b_\mm := \sum_{\kk \subseteq \mm } (-1)^{|\kk | } \binom{\mm }{\kk }_\theta  a_{\kk}
.\end{align*}
Then, their formal generating functions $ G_A(\zzz ):=\sum_{\kk\in \PN} a_{\kk }\Psi_{\kk }^\theta (\zzz )$ and $ G_B(\zzz ):=\sum_{\kk\in \PN} b_{\kk }\Psi_{\kk }^\theta (\zzz ) $ satisfy the relation 
\begin{align*}
  e^{-| \zzz |}G_A(\zzz )=G_B(- \zzz )
.\end{align*}
\begin{proof}
  This lemma immediately follows from \eqref{:25b}.
\end{proof}

\end{lemma}

\begin{proposition}\label{p:315}
For all $\mm \in \PN$ and $\xxx \in \R^N$, we have 
\begin{align}\label{:316a}
P_{\mm}^{\theta, \ip } (\xxx +\theta \delta )=(a)_{\mm}^\theta  P_{\mm}^\theta (\ctil \one  ) \sum_{\kk \subseteq \mm } (-1)^{|\kk |}\binom{\mm }{\kk }_\theta  M_\kk^{\theta ,a, c}(\xxx )
.\end{align}
\begin{proof}
It suffices to prove \eqref{:316a} for $\xx \in \PN$.
Applying \lref{l:314}, and using \eqref{:27a}, we have 
\begin{align*}
\sum_{\mm \in \PN } \Big( \sum_{\kk \in \PN } (-1)^{|\kk |}\binom{\mm}{\kk }_\theta  M_{\kk }^{\theta,a,c} (\xx ) \Big) \Psi_{\mm }^{\theta}(- \zzz )
&
=e^{-|\zzz | } \sum_{\kk \in \PN } M_{\kk }^{\theta, a,c} (\xx )\Psi_{\kk}^\theta (\zzz )
\\&
=\frac{1}{(a)_\xx^\theta \Psi_{\xx}^\theta (\one  ) } \LL_{\xx}^{\theta, a} \Big(\frac{\zzz }{\ctil } \Big)
\\&
=\sum _{\mm \subseteq \xx } (-1)^{|\mm |} \binom{\xx}{\mm } _\theta \frac{1}{(a)_\mm^\theta } \Phi_{\mm }^\theta \Big(\frac{\zzz }{\ctil } \Big)
.\end{align*}
Therefore, comparing the coefficients of $\Phi_{\mm}^\theta (\zzz )$, we complete the proof.
\end{proof}
\end{proposition}

\begin{remark}
From \eqref{:316a}, we have 
  \begin{align*}&
\sum_{\xx \in\PN }  P_{\mm}^{\theta, \ip } (\xx + \theta \delta ) e^{-\ctil |\zzz | } \Psi_{\xx }^\theta(\ctil \zzz )
\\&\quad 
=
(a)_{\mm}^\theta  P_{\mm}^\theta (\ctil \one  ) \sum_{\kk \subseteq \mm } (-1)^{|\kk |}\binom{\mm }{\kk }_\theta  
\sum _{\xx \in \PN } M_{\kk}^{\theta, a,c} (\xx  ) e^{-\ctil |\zzz | } \Psi_{\xx}^\theta (\ctil \zzz )    
. \end{align*}
Evaluating the left-hand side by \eqref{:25b}, and the right-hand side by \eqref{:27a} together with the duality \eqref{:31e}, we obtain the inversion formula for the Laguerre polynomials
\begin{align}\label{:316b}
P_{\mm}^{\theta } (\zzz  )
=(a)_{\mm}^\theta P_{\mm}^\theta ( \one  ) \sum_{\kk \subseteq \mm } (-1)^{|\kk |}\binom{\mm }{\kk }_\theta \frac{1}{(a)_{\kk }^\theta \Psi_{\kk}^\theta (\one )} \LL_\kk^{\theta,a } (\zzz  )
,\end{align}
which was proved in \cite[(4.14b)]{BaF97}.
By reversing the above computation and comparing coefficients at the final step, we obtain \eqref{:316a} from \eqref{:316b}.
\end{remark}

\section{Generalities on Markov kernels and Markov processes}\label{s:4}
\subsection{Markov kernels}\label{s:41}
Here we recall basic concepts of kernels and semigroups (see, for instance, \cite[Section 2]{BoO12} for a more comprehensive discussion).
Let $E$ and $E'$ be Borel subsets of Euclidean spaces.
Consider a function $\Lambda(x,A)$, where $x\in E $ and $A$ is a Borel subset of $E'$.
We say that $\Lambda $ is a Markov kernel from $E$ to $E'$, denoted by $\Lambda : E \dashrightarrow E' $, if the following two conditions hold:
\begin{itemize}
\item
 $\Lambda (x ,\cdot )$ is a Borel probability measure on $E'$ for any $x\in E$.

\item
$\Lambda(\cdot, A) $ is a Borel measurable function on $E$ for any Borel subset $A\subset E'$.
\end{itemize}

Let $\mathcal{B}_b(E) $ and $\mathcal{B}_b(E') $ denote the Banach spaces of $\mathbb{R} $-valued bounded Borel measurable functions on $E$ and $E'$, respectively, equipped with the supremum norm.
Let $C_\infty(E)$ be the space of all continuous functions on $E$ vanishing at infinity.
A Markov kernel $\Lambda : E \dashrightarrow E' $ defines a linear operator $\mathcal{B}_b(E') \to \mathcal{B}_b(E)$ by $(\Lambda f)(x)=\int_{E'}\Lambda (x, dy) f(y)$.
The Markov kernel $\Lambda  $ is called Feller if the induced operator satisfies $\Lambda(C_\infty(E')) \subset C_\infty(E) $.

Let $\mathcal{M}_p(E)$ and $\mathcal{M}_p(E')$ denote the sets of all probability measures on $E$ and $E'$, respectively.
For $\mu \in \mathcal{M}_p(E)$, define a probability measure $\mu \Lambda  \in \mathcal{M}_p(E')$ by $\mu \Lambda (A)=\int_{E}\Lambda (x, A ) \mu (dx) $.
Thus, a Markov kernel $\Lambda : E \dashrightarrow E' $ induces a map $\mathcal{M}_p(E) \to \mathcal{M}_p(E')$.

A Markov semigroup on $E$ is defined as a family of Markov kernels $\{ T_t\}_{t \ge 0}$ on $E$ such that $T_0 =I$ and $T_t T_s=T_{t+s}$.
This semigroup induces a semigroup on $\mathcal{B}_b(E)$.
A Markov semigroup on $E$ is said to be Feller--Dynkin if the following two conditions hold:

\begin{itemize}
\item 
The induced map $T_t: \mathcal{B}_b(E) \to \mathcal{B}_b(E)$ satisfies $T_t ( C_\infty(E) ) \subset C_\infty(E) $ for any $t\ge 0$.

\item
The induced semigroup on $C_\infty(E)$ is strongly continuous, namely,  $\lim_{t\to 0} T_t f =f $ in $C_\infty(E)$ holds for any $f\in C_\infty(E)$.
\end{itemize}

We remark that Feller--Dynkin semigroups are often called simply Feller.
To avoid confusion, we distinguish between the Feller and Feller--Dynkin properties.
In our terminology, the Feller property means just the property of mapping $C_\infty(E)$ into itself;  the Feller--Dynkin property
additionally requires strong continuity on $C_\infty(E)$.

\subsection{Continuous-time Markov chains} \label{s:42}
We recall a construction theorem for continuous-time Markov chains on a countable set $E$ equipped with the discrete topology.
A matrix $Q:=[Q(x,y)]_{x,y \in E}$ is called a Q-matrix if the following conditions hold:
\begin{alignat*}{2}
  Q(x,y) &\ge 0
    &\qquad& \text{for all } x, y \in E \text{ with }x\neq y, \\
 -Q(x,x) &= \sum_{y;\,y\neq x} Q(x,y)<\infty
    &\qquad& \text{for all } x \in E .
\end{alignat*}

For a Q-matrix $Q$, there exists a substochastic semigroup $\{ P_{\mathrm{min} ,t }\}_{t\ge 0}$ which is the minimal solution to Kolmogorov's backward and forward equations 
\begin{align*}
  \OD{}{t} P_t =Q P_t , \qquad   \OD{}{t} P_t = P_t Q
\end{align*}
with $P_0=I$.
We say that a Q-matrix is regular if the corresponding minimal solution  $\{ P_{\mathrm{min} ,t }\}_{t\ge 0}$ is stochastic.
If a Q-matrix is regular, then the semigroup $\{ P_{\mathrm{min} ,t }\}_{t\ge 0}$ is the unique stochastic solution to both Kolmogorov's backward and forward equations with initial condition $P_0=I$.

The following version of the construction theorem suffices for our purpose, as in \cite[Theorem 2.3]{BoO13}; see \cite[Chapter 8, Theorem 3.1]{EtK86} for a more general setting.
\begin{lemma}\label{l:41}
Assume that a Q-matrix $Q$ has finitely many non-zero entries in every row and every column.  
Assume that there exist a constant $C>0$ and strictly positive functions $\gamma , \eta  $ on $E$, both tending to infinity at infinity, such that, for any $x \in E $,  
\begin{align}\label{:41s}
-Q(x,x) \le C \gamma (x ) 
\\\label{:41t}
\Big( Q \frac{1}{\gamma } \Big) (x) \le \frac{C}{\gamma (x)}
\\\label{:41u}
(Q \eta ) (x) \le C \eta(x)  
,\end{align}
where $Qf (x):=\sum_{y \in E} Q(x,y)f(y)=\sum_{y\in E, y \neq x}Q(x,y) (f(y)-f(x))$.
Then, $Q$ is regular and therefore determines a Markov semigroup $\{ T_t \}_{t\ge 0}$.
Furthermore, $\{ T_t \}_{t\ge 0}$ is Feller--Dynkin.
\end{lemma}

\section{Gateways between Laguerre and Meixner processes}\label{s:5}
In this section, we work with the parameters $\alpha,\sigma>0$ in place of $a,c$, and apply the formulae of \sref{s:3} after the reparametrisation
\begin{align}\label{:41f}
a = \theta(\alpha +N-1)
, \qquad 
c=\frac{\sigma}{1+\sigma}
\end{align}
together with the coordinate change given by the bijection $\map{\varphi }{\nWNd }{ \PN}$  defined by 
\begin{align}\label{:41e}
\varphi ( \xx ) = \xx^R - \delta
,\end{align}
where $\xx^R =(x_N, \ldots , x_1)$ for $\xx= (x_1,\ldots, x_N)$.

\subsection{Meixner processes}
Let $Q_{\theta, \alpha ,\sigma}^N := [Q_{\theta, \alpha ,\sigma}^N (\xx,\yy ) ]_{\xx,\yy \in \nWNd }$ be the matrix associated with the generator $\mathbf{A}_{\theta,\alpha ,\sigma }^{N}$ defined by \eqref{:11g}, that is,
\begin{align}\label{:41k}
  Q_{\theta, \alpha , \sigma }^{N} (\xx, \xx+ \eee _i )&:= \sigma q_{\theta, \alpha }^{N} (\xx, \xx+ \eee _i )
,\\ \label{:41l}
  Q_{\theta, \alpha , \sigma }^{N} (\xx, \xx- \eee_i )&:= (1+\sigma) q_{\theta, \alpha }^{N} (\xx, \xx- \eee_i ) 
,\\  \notag
 - Q_{\theta, \alpha , \sigma}^N (\xx,\xx  )
 &:= \sum_{i=1}^{N } ( Q_{\theta, \alpha , \sigma }^{N} (\xx, \xx+ \eee_i )+Q_{\theta, \alpha , \sigma }^{N} (\xx, \xx- \eee_i ))
 \\& \label{:41m}
 =  \frac{1+2\sigma}{1+\sigma }\Big(|\xx |-\frac{N(N-1)}{2 } \Big)+\frac{\sigma}{1+\sigma } \theta(\alpha + N-1) N 
,\end{align}
and we set $Q_{\theta, \alpha ,\sigma}^N (\xx,\yy  ):=0$ otherwise.
Here, the last equality in \eqref{:41m} follows since \eqref{:33e} and \eqref{:33f} imply  
\begin{align}\label{:52c}
\sum_{i=1}^{N } Q_{\theta, \alpha , \sigma }^{N} (\xx, \xx+ \eee_i )&=\frac{\sigma}{1+\sigma }\Big(|\xx |-\frac{N(N-1)}{2 } + \theta(\alpha + N-1) N \Big)
, \\\label{:52d}
  \sum_{i=1}^{N } Q_{\theta, \alpha , \sigma }^{N} (\xx, \xx- \eee_i ) &=|\xx| -\frac{N(N-1)}{2}
.\end{align}

\begin{lemma}\label{l:51}
Assume that $\theta, \alpha , \sigma >0$.
Then, the matrix $Q_{\theta, \alpha ,\sigma}^N $ is a regular Q-matrix.
In particular, there exists a Markov semigroup $\{ \mathbf{T}_{\theta, \alpha , \sigma , t}^{N}\}_{t\ge 0}$ on $\nWNd $ associated with $Q_{\theta, \alpha ,\sigma}^N$.
Furthermore, the semigroup $\{ \mathbf{T}_{\theta, \alpha , \sigma , t}^{N}\}_{t\ge 0}$ is Feller--Dynkin.
\begin{proof}
  To see the positivity of $ Q_{\theta, \alpha , \sigma }^{N} (\xx, \xx\pm  \eee_i )$, from 
\begin{align*}
  \frac{x_i-x_j-i+j-\theta(j-i\pm 1)}{x_i-x_j-i+j-\theta(j-i) }=1\mp  \frac{\theta }{x_i-x_j-i+j-\theta(j-i) }
,\end{align*}
it suffices to show that 
\begin{align}\label{:41i}
  |x_i-x_j-i+j-\theta(j-i) |\ge \theta 
.\end{align}
If $i <j$, then we have $x_i-x_j-i+j \le 0$ and $-\theta(j-i)\le -\theta$.
If $i > j$, then we have $x_i-x_j-i+j \ge 0$ and $-\theta(j-i)\ge \theta$.
Thus, we have \eqref{:41i}, which proves that $Q_{\theta, \alpha ,\sigma}^N $ is a Q-matrix.
From \eqref{:52c} and \eqref{:52d}, it is easy to see that the conditions \eqref{:41s}--\eqref{:41u} hold for $\gamma (\xx )=\eta (\xx ) := |\xx|+1$.
Therefore, by \lref{l:41}, we finish the proof.
\end{proof}
\end{lemma}

Let $\XX_t=(X_t^1, \ldots, X_t^N )$ be the Meixner process associated with $Q_{\theta, \alpha ,\sigma}^N $.
The distribution of the Meixner process at time $t$ has exponential moments.

\begin{lemma}\label{l:53}
Let $t\ge 0$ and $\xx \in \nWNd $.
Then, $\E_{\xx } [e^{\ep |\XX _t |}]<\infty $ holds for sufficiently small $\ep>0$.
\begin{proof}
There exist constants $\lambda, \nu >0$ such that $Q_{\theta, \alpha, \sigma}^N (\xx, \xx + \eee_N) \le \lambda x_N + \nu $.
Let $\{ Z_t \}_{t\ge 0}$ be the pure birth process on $\Z _{\ge }$ with generator $\mathbf{A} := (\lambda x + \nu )\nabla_{+} $. 
We set $u(t,x)=\E_{x} [e^{\ep Z_t}]$, which satisfies the backward equation
\begin{align*}
 \partial_{t} u(t,x)=[\mathbf{A}  u(t,\cdot)] (x), \quad u(0,x)=e^{\ep x}
.\end{align*}
Define $H(s, \xx )= u(t-s, x_N)$ for $0\le s\le t $.
Since $H$ depends only on the last coordinate and $u (t-s, x)$ is increasing in $x$, we have
\begin{align*}
[ Q_{\theta, \alpha, \sigma}^N H(s, \cdot ) ](\xx )
 &\le Q_{\theta, \alpha, \sigma}^N  (\xx, \xx + \eee_N)( H(s, \xx +\eee_N ) -  H(s, \xx ))
 \\&
 \le [\mathbf{A}  u(t-s, \cdot  ) ](x_N)
 \\&
 =
 -\partial_{s} H(s, \xx )
.\end{align*}
Therefore, for a stopping time $\tau_m:= \inf \{t\ge 0 \,;\, X_t^N \ge m\} $, Dynkin's formula yields
\begin{align*}
  \E [ H(t\wedge \tau_m, \XX_{t \wedge \tau_m }) ] 
  =H(0, \XX_0) +\E \Big[\int _0^{t\wedge \tau_m}  [(  \partial_{s} + Q_{\theta, \alpha, \sigma}^N) H(s, \cdot )] ( \XX_{s}) ds \Big]
\le H(0, \XX_0)
.\end{align*}
Hence, letting $m\to \infty $ with Fatou's lemma, we have  $  \E_{\xx } [e^{\ep X_{t}^N}] \le u(t, x_N)$ using  $\tau_m \uparrow \infty $.
Thus, the exponential moments of the $N$-th particle $X_t^{N} $ are bounded by those of $Z_t$.

Since the probability generating function for $Z_t$ is given by
\begin{align*}
  \E_x [s ^{Z_t }] = \frac{s^x }{(s+e^{\lambda t}(1-s ))^{x +\frac{\nu}{\lambda }}}
\end{align*}
(see \cite[Section 3.2]{AndW91} for example),  $\E_x [e ^{\ep Z_t }]<\infty $ for sufficiently small $\ep$.
Combining this with $\E_{\xx } [e^{\ep |\XX _t |}] \le \E_{\xx } [e^{\ep N X _t ^N }] $, we complete the proof.
\end{proof}
\end{lemma}

Under the parameter change \eqref{:41f} and the coordinate change \eqref{:41e}, the normalised orthogonality measure in \pref{p:313} induces a probability measure on $\nWNd$, which we call the Meixner ensemble and denoted by $\mathbf{m}_{\theta, \alpha , \sigma }^N$.
Explicitly, for $\xx\in\nWNd$, we have
\begin{align}
\mathbf{m}_{\theta, \alpha ,\sigma }^N (\xx )= (\theta(\alpha +N-1) )_{\varphi (\xx) }^\theta \Psi_{\varphi (\xx ) } ^\theta \Big( \Big(\frac{\sigma }{1+\sigma } \Big) \one  \Big) \Big(1-\frac{\sigma }{1+\sigma } \Big)^{N\theta(\alpha +N-1) }
.\end{align}

\begin{remark}
For $\theta=1/2, 1,2$, we have the following expression, where $y_i=2x_i-i+1$ and $z_i=x_i+i-1$:
\begin{align*}
 \mathbf{m}_{1/2, \alpha ,\sigma }^N (\xx )
 &
 \propto \prod_{1\le i <j\le N } (y_j -y_i )\prod_{j=1}^{N} \frac{\Gamma (\frac{y_j+\alpha }{2} ) }{\Gamma (\frac{y_j+2 }{2} ) }  \Big(\frac{\sigma }{1+\sigma } \Big) ^{y_j /2 }
, \\
 \mathbf{m}_{1, \alpha ,\sigma }^N (\xx )
 &
 \propto \prod_{1\le i <j\le N } (x_j -x_i )^2 \prod_{j=1}^{N} \frac{(\alpha )_{x_j}}{x_j !}  \Big(\frac{\sigma }{1+\sigma } \Big) ^{x_j}
,\\
 \mathbf{m}_{2, \alpha ,\sigma }^N (\xx )
 &
 \propto \prod_{1\le i <j\le N } (z_j -z_i )^2(z_j -z_i -1)(z_j-z_i+1 ) \prod_{j=1}^{N} \frac{(2 \alpha )_{z_j}}{z_j !}  \Big(\frac{\sigma }{1+\sigma } \Big) ^{z_j}
.\end{align*}
The measure $ \mathbf{m}_{1, \alpha ,\sigma }^N $ corresponds to the Meixner ensemble introduced in \cite{Joh01}, and $ \mathbf{m}_{2, \alpha ,\sigma }^N $ corresponds to the Meixner symplectic ensemble introduced in \cite[Proposition 4.1]{BoS09} under the coordinate change $\mathbf{z}=\xx+\delta^R$.
If $N$ is even, $ \mathbf{m}_{1/2, \alpha ,\sigma }^N $ corresponds to the Meixner orthogonal ensemble  introduced in \cite[Proposition 4.2]{BoS09} under the coordinate change $\mathbf{y}=2\xx - \delta^R$.
\end{remark}

\begin{proposition}\label{p:53}
The Meixner ensemble $\mathbf{m}_{\theta, \alpha ,\sigma  }^N$ is the unique invariant distribution for the Meixner process $\{ \mathbf{T}_{\theta, \alpha , \sigma , t}^{N}\}_{t\ge 0}$.
Moreover, it is a symmetrising measure for this process.
\begin{proof}
Since the irreducibility of the process follows immediately from the transition rates defined by \eqref{:11g}, the invariant distribution is unique.
Then, it remains to prove the following detailed balance condition: for $\xx, \yy \in\nWNd $, 
\begin{align}\label{:43a}
\mathbf{m}_{\theta, \alpha , \sigma  }^N (\xx )   Q_{\theta, \alpha ,\sigma}^N(\xx, \yy )  = \mathbf{m}_{\theta , \alpha ,\sigma  }^N (\yy )Q_{\theta, \alpha ,\sigma}^N(\yy, \xx )
. \end{align}
It suffices to establish \eqref{:43a} for $\yy=\xx+\eee_i$, for all $\xx\in \nWNd$, and this identity follows from \eqref{:315a}.
Thus, we complete the proof.
\end{proof}
\end{proposition}

\subsection{Properties of Markov kernels}

For $\theta ,\alpha, \sigma>0$, we define kernels $ \Lambda_{\theta, \sigma }^{N} : \nWN \dashrightarrow \nWNd $ and $\Lambda_{\theta, \alpha , \sigma}^{*, N} : \nWNd \dashrightarrow \nWN $ by
\begin{align}
  \label{:41c}
\Lambda_{\theta, \sigma }^{N} (\xxx, \yy ) &= \Psi_{\varphi (\yy ) }^\theta (\sigma \xxx ) e^{- \sigma  |\xxx | }
,\\
\label{:41d}
\Lambda_{\theta, \alpha , \sigma}^{*, N} (\xx, d\yyy ) &=\frac{1}{Z_{\theta, \alpha, \sigma } }\frac{(1+\sigma)^{|\xx| } \Phi_{\varphi (\xx ) }^{\theta} (\yyy ) }{ (\theta(\alpha +N-1 ))_{\varphi (\xx ) }^\theta }  e^{-(1+\sigma)|\yyy | } \prod_{j=1}^N y_j^{\theta \alpha -1 }  \Delta(\yyy )^{2\theta }  d\yyy 
.\end{align}
Here, the normalising constant is given by
\begin{align*}
\frac{1}{Z_{\theta, \alpha ,\sigma } } =\frac{N! }{\mathcal{N}_{\theta, \theta(\alpha +N-1)}} (1+\sigma)^{\theta N (\alpha + N-1) -N(N-1)/2} 
.\end{align*}
The kernels $ \Lambda_{\theta, \sigma }^{N} $ and $\Lambda_{\theta, \alpha , \sigma}^{*, N}  $ are Markov since their total masses are $1$ by taking $\mm=0$ in \eqref{:62a} and \eqref{:52a} below, and positivity follows from \lref{l:1}.

\begin{remark}
The distribution $\Lambda_{\theta, \sigma }^N (\xxx , \cdot ) $ is a multivariate analogue of the Poisson distribution.
Furthermore, by rewriting \eqref{:25b}, we have 
\begin{align}\label{:56a}
\sum_{\yy \in \nWNd }P_{\kk}^{\theta, \ip }(\varphi (\yy )  +\theta \delta) e^{-\sigma | \xxx |} \Psi_{\varphi (\yy ) }^\theta(\sigma \xxx ) =P_{\kk } ^\theta (\sigma \xxx )
,\end{align}
which is a multivariate analogue of the formula stating that the $k$-th falling factorial moment of the Poisson distribution with parameter $\sigma$ is given by the monomial $\sigma^k $.
\end{remark}

Let $\tilde{M}_\mm^{\theta, \alpha , \sigma }$ be the function on $\nWNd $ defined by
\begin{align}
\tilde{M}_\mm^{\theta, \alpha , \sigma } (\xx )=  M_\mm^\theta \Big(\varphi(\xx ) ;\, \theta(\alpha+N-1),\frac{\sigma}{1+\sigma } \Big)
.\end{align}
Then, from \eqref{:27a} and \eqref{:36d}, for all $\mm\in \PN $, $\xxx\in \nWN $, and $\xx \in \nWNd $, we have
\begin{align}\label{:62a}
[\Lambda_{\theta, \sigma }^{N} \tilde{M}_{\mm}^{\theta, \alpha, \sigma } ](\xxx )
&= \frac{1}{(\theta(\alpha +N-1))_{\mm }^\theta \Psi_{\mm}^{\theta}( \one ) } \LL_{\mm}^{\theta, \theta(\alpha +N-1)}(\xxx ) 
,\\ \label{:52a}
 [\Lambda_{\theta,\alpha , \sigma }^{*, N}  \LL_{\mm}^{\theta,\theta(\alpha +N-1)} ] (\xx )
&=\Big( \frac{\sigma }{1+\sigma } \Big)^{|\mm |}  (\theta(\alpha +N-1) )_{\mm }^\theta \Psi_{\mm}^{\theta}( \one  )  \tilde{M}_{\mm } ^{\theta , \alpha ,\sigma} ( \xx )
.\end{align}
Furthermore, \eqref{:33c} yields
\begin{align}\label{:41h}
  \mathbf{A}_{\theta, \alpha ,\sigma }^{N} \tilde{M}_\mm^{\theta, \alpha ,\sigma } (\xx )= -|\mm | \tilde{M}_\mm^{\theta, \alpha ,\sigma} (\xx )
.\end{align}
We prepare moment estimates for the Markov kernels to apply the moment problem.
\begin{lemma}
Assume that $\theta>0$.
Then, for every $\ep >0$, there exist $\eta >0 $ and $\Ct{c;1}>0$ such that 
 \begin{align} \label{:64a}
[ \Lambda_{\theta, \sigma }^{N} e^{\eta |\cdot | } ](\xxx ) 
&
\le \cref{c;1} e^{\ep |\xxx |} \qquad \text{ for all } \xxx \in \nWN
,\\  \label{:64b}
[ \Lambda_{\theta, \alpha , \sigma }^{*,N} e^{\eta |\cdot | } ](\xx )
&
\le \cref{c;1} e^{\ep |\xx |} \qquad  \text{ for all } \xx \in \nWNd
. \end{align}
\begin{proof}
By \eqref{:25b} and  \eqref{:34d}, we have  
\begin{align*}
[ \Lambda_{\theta, \sigma }^{N} e^{\eta |\cdot | } ](\xxx ) 
&=e^{- \sigma  |\xxx |  + \eta \frac{N(N-1)}{2}}  e^{\sigma e^{\eta } |\xxx | }
,\\
[ \Lambda_{\theta, \alpha , \sigma }^{*,N} e^{\eta |\cdot | } ](\xx )&= \Big(\frac{1+\sigma}{1+\sigma- \eta } \Big)^{N\theta(\alpha+N-1)+|\xx | -N(N-1)/2}
,\end{align*}
which yields \eqref{:64a} and \eqref{:64b}, respectively.

\end{proof} 
\end{lemma}

\begin{proposition}
For $\theta>0$, the Markov kernels $\Lambda_{\theta, \sigma }^{N}  $ and $\Lambda_{\theta, \alpha , \sigma }^{*,N}  $ are Feller.
  Namely, we have
  \begin{align}\label{:65a}
    \Lambda_{\theta, \sigma }^{N} f &\in C_{\infty }(\nWN ) \text{ for any } f \in C_{\infty }(\nWNd )
  ,  \\ \label{:65b}
   \Lambda_{\theta, \alpha , \sigma }^{*,N}  f  &\in C_{\infty }(\nWNd ) \text{ for any } f \in C_{\infty }(\nWN )
 . \end{align} 

 \begin{proof}
Note that the Feller property \eqref{:65b} simply means that the kernel preserves functions vanishing at infinity, which follows readily from \eqref{B:1a} by the dominated convergence theorem.

For $ f\in C_{\infty }(\nWNd )$, we choose finite sets $F_n$ such that $\limi{n}\|f-f_n\|=0$ for $f_n=f \mathbf{1}_{F_n}$.
Then we obtain
\begin{align}\label{:65c}
\| \Lambda_{\theta, \sigma }^N f - \Lambda_{\theta, \sigma }^N f_n \|
\le
\sup_{\xxx \in\nWN } \sum_{\yy \in \nWNd } \Lambda_{\theta, \sigma }^N(\xxx, \yy)| f(\yy) - f_n(\yy)|
\le \| f -f_n\| \xrightarrow[n\to\infty]{}  0
.\end{align}
Since we have $\Lambda_{\theta, \sigma }^N f_n \in C_{\infty }(\nWN )$ by the definition of $\Lambda_{\theta, \sigma }^N$, we obtain \eqref{:65a} from \eqref{:65c}.
\end{proof}
\end{proposition}

\begin{lemma}\cite{deJ03,DLY02}\label{l:59}
For a probability measure  $\mu $ on $\R_{\ge}^N$, assume that $\int_{\R_{\ge}^N} e^{\ep |\xxx |} d \mu (\xxx ) <\infty $ for some $\ep >0$.
Then, the measure $\mu$ is uniquely determined by its moments.
\end{lemma}


We have the following lemma, a multivariable version of  \cite[Corollary 3.1.6]{BoO13b}.

\begin{lemma}\label{l:66}
The function space $\Lambda_{\theta, \sigma }^{N} (C_{\infty }(\nWNd ) )$ is dense in $C_{\infty }(\nWN )$.
  \begin{proof}

We first note that, for  $\mu \in \mathcal{M}_{p}(\nWN )$, the finite measure $e^{- \sigma |\xxx |} d \mu (\xxx ) $ is uniquely determined by its moments from \lref{l:59} and a trivial identity $\int_{\nWN} e^{\sigma |\xxx |} e^{- \sigma |\xxx |} d \mu (\xxx ) =1$.

To complete the proof of this lemma, it suffices to show that the linear span of $\{ \Lambda_{\theta, \sigma }^{N}(\cdot ,\yy ) \}_{\yy \in \nWNd }$ is dense in $C_{\infty }(\nWN )$.
This is reduced to proving that if a finite signed measure $\nu$ on $\nWN$ satisfies  
\begin{align*}
\int_{\nWN } \nu(d\xxx) \Psi_{\varphi (\yy ) }^\theta (\sigma \xxx) e^{-\sigma | \xxx |} =0
,\end{align*}
for all $\yy \in \nWNd $, then $\nu =0$. 
We set $\nu= \nu_1 -\nu_2$ for finite positive measures $\nu_1,\nu_2$.
Then, the assumption on $\nu $ implies that $e^{- \sigma |\xxx |} d \nu_1 (\xxx )  $ and $e^{- \sigma |\xxx |} d \nu_2 (\xxx )  $ have the same symmetric polynomial moments.
Symmetrising these finite positive measures $e^{- \sigma |\xxx |} d \nu_i (\xxx )  $ ($i=1,2$) and applying \lref{l:59}, we conclude that they coincide, since both original measures are supported on the Weyl chamber.
Hence, we have $\nu_1 =\nu_2$, which yields $\nu=0$.
Thus, we finish the proof.
\end{proof}
\end{lemma}

\subsection{Proof of \tref{t:11}}
Following the approach of \cite{RaS18}, we prove \tref{t:11} by means of the moment problem.
However, the argument is simpler in our setting than in \cite{RaS18}, since the families of multivariable Meixner and Laguerre polynomials diagonalise the generators of the processes under consideration.

Let $\{\mathsf{X} _t\}_{t\ge0}$ denote the Laguerre process associated with the semigroup $\{T_{\theta,\alpha,t}^{N}\}_{t\ge0}$.
The Laguerre process has the exponential moments; see \cite[Proof of Proposition B.2]{KaS26}.
\begin{lemma}\label{l:67}
Let $t\ge 0$ and $\xxx \in \nWN $.
Then, for sufficiently small $\ep>0$, we have $\E_\xxx [e^{\ep |\mathsf{X} _t |}]<\infty $.
\end{lemma}

For $m \in \mathcal{M}_p (\nWN) $, define a symmetrised probability measure $m^{\sym }$ on $\R _{\ge }^N $ by
\begin{align*}
  m^{\sym } (dz_1,\ldots , dz_{N})= \frac{1}{N!} m (dz_{(1)},\ldots ,dz_{(N)})
,\end{align*}
where $z_{(1)}\le z_{(2)} \le \cdots \le z_{(N)}$ are the order statistics of $(z_1, z_2,\ldots, z_N)$.
Furthermore, for a (not necessarily symmetric) function $f$ on $ \R_{\ge}^N $, we have
\begin{align}\label{:42b}
  \int_{\R_{\ge}^N } f(\zzz ) m^\sym( d\zzz ) =\int_{\R_{\ge}^N} f^\sym (\zzz ) m^\sym(d \zzz ) =  \int_{\nWN } f^\sym (\zzz ) m(d \zzz )
,\end{align}
where $f^\sym $  is a symmetric function given by $f^\sym (\zzz )=\frac{1}{N!} \sum_{\sigma\in \mathfrak{S}_N}f(z_{\sigma(1) },\ldots, z_{\sigma(N) })$.

\medskip
\noindent
{\em Proof of \tref{t:11} } \quad 
By It\^o’s formula for $\LL_{\mm }^{\theta, \theta(\alpha+N-1)} (\mathsf{X}_t)$ and then taking expectations, justified by \lref{l:67}, and using \eqref{:31a}, we have
\begin{align*}
[\mathsf{T}_{\theta, \alpha, t}^{N} \LL_{\mm }^{\theta, \theta(\alpha+N-1)} ](\xxx ) 
&= \LL_{\mm }^{\theta, \theta(\alpha+N-1)} (\xxx )  + \int_0^t [\mathsf{T}_{\theta, \alpha, s}^{N} \mathsf{A}_{\theta, \theta(\alpha+N-1) }^{N} \LL_{\mm }^{\theta, \theta(\alpha+N-1)} ] (\xxx ) ds
 \\& \notag
 = \LL_{\mm }^{\theta, \theta(\alpha+N-1)} (\xxx )  -|\mm| \int_0^t [\mathsf{T}_{\theta, \alpha, s}^{N} \LL_{\mm }^{\theta, \theta(\alpha+N-1)} ] (\xxx ) ds
,\end{align*}  
which yields 
\begin{align}\label{:69a}
[\mathsf{T}_{\theta, \alpha, t}^{N} \LL_{\mm }^{\theta, \theta(\alpha+N-1)} ](\xxx ) 
&= e^{-|\mm |t } \LL_{\mm }^{\theta, \theta(\alpha+N-1)} (\xxx ) 
.\end{align} 
Similarly, by Dynkin's formula for $\tilde{M}_{\mm }^{\theta, \alpha, \sigma }(\XX _t)$ with the moment estimates in \lref{l:53} and using \eqref{:41h}, we have
\begin{align*}
[\mathbf{T}_{\theta, \alpha,\sigma, t}^{N} \tilde{M}_{\mm }^{\theta, \alpha, \sigma } ](\xx ) 
&=  \tilde{M}_{\mm }^{\theta, \alpha, \sigma } (\xx )   + \int_0^t [\mathbf{T}_{\theta, \alpha, \sigma, s}^{N} \mathbf{A}_{\theta, \alpha, \sigma  }^{N}  \tilde{M}_{\mm }^{\theta, \alpha, \sigma } ](\xx )  ds
 \\& \notag
 = \tilde{M}_{\mm }^{\theta, \alpha, \sigma }(\xx) -|\mm| \int_0^t  [\mathbf{T}_{\theta, \alpha, \sigma, s}^{N}  \tilde{M}_{\mm }^{\theta, \alpha, \sigma } ](\xx ) ds
,\end{align*}  
which yields 
\begin{align}\label{:69b}
[\mathbf{T}_{\theta, \alpha, \sigma, t}^{N} \tilde{M}_{\mm }^{\theta, \alpha, \sigma } ](\xx ) &= e^{-|\mm |t } \tilde{M}_{\mm }^{\theta, \alpha, \sigma } (\xx )
.\end{align} 

Combining \eqref{:69a} and \eqref{:69b} with \eqref{:62a}, we obtain 
\begin{align}\notag 
 [\mathsf{T}_{\theta, \alpha, t}^{N}  \Lambda_{\theta, \sigma}^{N} \tilde{M}_{\mm } ^{\theta , \alpha ,\sigma} ] ( \xxx )
 &=\frac{e^{-|\mm |t }}{(\theta (\alpha +N-1))_{\mm }^\theta \Psi_{\mm}^{\theta}( \one  ) } \LL_{\mm }^{\theta, \theta(\alpha+N-1)} (\xxx ) 
\\& \label{:69c}
=[ \Lambda_{\theta, \sigma}^{N} \mathbf{T} _{\theta, \alpha, \sigma, t}^{N}  \tilde{M}_{\mm } ^{\theta , \alpha ,\sigma} ](\xxx )
.\end{align}
Similarly, combining \eqref{:69a} and \eqref{:69b} with \eqref{:52a}, we obtain 
\begin{align}\label{:69d}
 [\mathbf{T} _{\theta, \alpha, \sigma, t}^{N} \Lambda_{\theta, \alpha, \sigma}^{N,*} \LL_{\mm }^{\theta, \theta(\alpha+N-1)} ](\xx )
 = [\Lambda_{\theta, \alpha, \sigma}^{N,*}  \mathsf{T}_{\theta, \alpha, t}^{N}  \LL_{\mm }^{\theta, \theta(\alpha+N-1)} ](\xx ) 
.\end{align}

We set $\mathsf{m}_1= \mathsf{T}_{\theta, \alpha, t}^{N}  \Lambda_{\theta, \sigma}^{N} (\xxx, \cdot)$ and $\mathsf{m}_2=\Lambda_{\theta, \sigma}^{N} \mathbf{T} _{\theta, \alpha, \sigma, t}^{N} (\xxx, \cdot) $, and $f(\xx )=\xx+(\theta-1)\delta^R$.
Then, from \eqref{:69c} and \pref{p:32}, we have $\mathbb{E}_{\mathsf{m}_1\circ f^{-1}}[p]=\mathbb{E}_{\mathsf{m}_2 \circ f^{-1}}[p]$ for every symmetric polynomial $p$. 
Hence, \eqref{:42b} yields that all moments of the symmetrised measures $(\mathsf{m}_1\circ f^{-1})^{\mathrm{sym}} $ and $(\mathsf{m}_2 \circ f^{-1})^{\mathrm{sym}}  $ are the same.
From this with \lref{l:59} and the exponential moment estimates \eqref{:64a}, \lref{l:53}, and \lref{l:67}, which remain valid after the fixed translation $f$, it follows that $(\mathsf{m}_1\circ f^{-1})^{\mathrm{sym}} =(\mathsf{m}_2 \circ f^{-1})^{\mathrm{sym}}   $.
Since both pushforward measures are supported on $\nWN$, equality of their symmetrisations implies $\mathsf{m}_1 \circ f^{-1}=\mathsf{m}_2 \circ f^{-1}$, and therefore $\mathsf{m}_1 =\mathsf{m}_2  $ by the injectivity of $f$.
Thus, this completes the proof of \eqref{:12a}.
By the same argument using \eqref{:69d} and the fact that the Laguerre polynomials form a basis for $\Q(\theta, a) [\xxx ]^{\mathfrak{S}_N}$, we obtain \eqref{:12b}.
Thus, we complete the proof of \tref{t:11}.
\qed
\medskip

\subsection{Applications of gateways}

Define the Laguerre ensemble on $\nWN $ by
\begin{align*}
  \mathsf{m}_{\theta, \alpha }^N (d\xxx )= \frac{1 }{\tilde{\mathcal{N}}_{\theta, \alpha }} e^{-|\xxx | } \prod_{k=1}^N x_k^{\theta \alpha -1}  \Delta(\xxx )^{2\theta }  d\xxx 
,\end{align*}
where
\begin{align*}
 \frac{1 }{\tilde{\mathcal{N}}_{\theta, \alpha }}=  \frac{N! }{\mathcal{N}_{\theta, \theta(\alpha+N-1)}}
.\end{align*}
It is known that $ \mathsf{m}_{\theta, \alpha }^N$ is an invariant and symmetrising measure for the Laguerre semigroup $ \{ \mathsf{T}_{\theta, \alpha ,t}^{N} \}_{t\ge 0}$.

Using the gateways, which transfer information from discrete processes to continuous processes, we obtain the following properties of the Laguerre semigroup.

\begin{proposition}\label{p:67}
Assume that $\theta \ge 1/2$ and $\alpha >0$.
Then the following hold.
  \begin{itemize}
\item[(i)]
The Laguerre ensemble  $  \mathsf{m}_{\theta, \alpha }^N$ is the unique invariant probability measure for the Laguerre semigroup $ \{ \mathsf{T}_{\theta, \alpha ,t}^{N} \}_{t\ge 0}$.

\item[(ii)]
The Laguerre semigroup $ \{ \mathsf{T}_{\theta, \alpha ,t}^{N} \}_{t\ge 0}$ is Feller--Dynkin.

\end{itemize}

\begin{proof}
Let $\mathsf{m}_i \in \mathcal{M}_p (\nWN)$ satisfy $\mathsf{m}_i =\mathsf{m}_i \mathsf{T}_{\theta, \alpha ,t}^{N}$ for $i=1,2$.
Then, the gateway \eqref{:12a} yields
\begin{align*}
 \mathsf{m}_i \Lambda_{\theta,\sigma}^N  
 =  \mathsf{m}_i \mathsf{T}_{\theta, \alpha ,t}^{N} \Lambda_{\theta,\sigma}^N 
 =  \mathsf{m}_i \Lambda_{\theta,\sigma}^N \mathbf{T}_{\theta,\alpha ,\sigma ,t}^N
.\end{align*}
Thus, $\mathsf{m}_i \Lambda_{\theta,\sigma}^N $ is an invariant distribution of $\{ \mathbf{T}_{\theta, \alpha ,\sigma ,t}^N \}_{t\ge 0}$.
Hence, the uniqueness of the invariant distribution for the Meixner process (\pref{p:53}) proves $\mathsf{m}_1 \Lambda_{\theta,\sigma}^N =\mathsf{m}_2 \Lambda_{\theta,\sigma}^N $.
Combining this with \lref{l:66}, we have $\mathsf{m}_1=\mathsf{m}_2$, which completes the proof of (i).

To show (ii), let $f\in C_{\infty }(\nWNd ) $.
Then, \lref{l:51} and \eqref{:65a} yield $\Lambda_{\theta,\sigma}^N \mathbf{T}_{\theta, \alpha ,\sigma ,t}^N f \in C_{\infty }(\nWN )$.
Combining this with \eqref{:12a}, we derive $ \mathsf{T}_{\theta, \alpha ,t}^{N} (\Lambda_{\theta,\sigma}^N f)\in C_{\infty }(\nWN )$. 
Therefore, from \lref{l:66}, we have $\mathsf{T}_{\theta, \alpha ,t}^{N} (C_{\infty }(\nWN ) ) \subset C_{\infty }(\nWN )$.
Furthermore, combining this with the continuity of paths of the process, the strong continuity follows from \cite[III Lemma 6.7]{RoW00}.
\end{proof}
\end{proposition}

The gateways yield also relations between invariant probability measures.
\begin{corollary}
For $\theta \ge 1/2$ and $\alpha ,\sigma >0$, the following identities of probability measures hold:
\begin{align}\label{:66a}
\mathsf{m}_{\theta, \alpha }^N \Lambda_{\theta,\sigma}^N  = \mathbf{m}_{\theta, \alpha  ,\sigma }^N 
,\\\label{:66b}
\mathbf{m}_{\theta, \alpha  ,\sigma }^N \Lambda_{\theta, \alpha, \sigma}^{N,*}  = \mathsf{m}_{\theta, \alpha }^N 
.\end{align}
\begin{proof}
Equation \eqref{:12a} implies that  
  \begin{align*}
\mathsf{m}_{\theta, \alpha }^N \Lambda_{\theta,\sigma}^N
=
\mathsf{m}_{\theta, \alpha }^N \mathsf{T}_{\theta, \alpha ,t}^{N} \Lambda_{\theta,\sigma}^N
=\mathsf{m}_{\theta, \alpha }^N \Lambda_{\theta,\sigma}^N\mathbf{T}_{\theta, \alpha ,\sigma ,t}^N
.\end{align*}
Hence, by the uniqueness of the invariant distributions established in \pref{p:53}, we obtain \eqref{:66a}.
Similarly, by \eqref{:12b}, we have $\mathbf{m}_{\theta, \alpha  ,\sigma }^N \Lambda_{\theta, \alpha, \sigma}^{N,*}  = \mathbf{m}_{\theta, \alpha ,\sigma }^N \Lambda_{\theta, \alpha, \sigma}^{N,*} \mathsf{T}_{\theta, \alpha, t}^{N} $, which, together with \pref{p:67} (i), yields \eqref{:66b}.

\end{proof}
\end{corollary}

The identities \eqref{:66a} and \eqref{:66b} can be written  explicitly as follows:
\begin{align}&\notag 
\frac{1 }{\tilde{\mathcal{N}}_{\theta, \alpha }}\int_{\nWN } \Psi_{\varphi(\yy ) }^\theta(\sigma \xxx )e^{-(1+\sigma)|\xxx | } \prod_{k=1}^N x_k^{\theta \alpha -1}  \Delta(\xxx )^{2\theta }   d\xxx 
\\&\label{:66c}
=(\theta(\alpha+N-1) )_{\varphi (\yy ) }^\theta \Psi_{\varphi (\yy ) } ^\theta \Big(\frac{\sigma }{1+\sigma } \one  \Big) \Big(1-\frac{\sigma }{1+\sigma } \Big)^{N\theta(\alpha+N-1) }
, \\ & \notag
\frac{1}{Z_{\theta, \alpha , \sigma }^N } \sum_{\xx \in \nWNd } \Psi_{ \varphi (\xx ) }^\theta \Big( \frac{\sigma }{1+\sigma }  \one  \Big) \Big(1-\frac{\sigma }{1+\sigma } \Big)^{Na }  (1+\sigma)^{|\xx | } 
\\& \notag \qquad\qquad\qquad   
\times\Phi_{\varphi (\xx )}^\theta(\yyy ) e^{-(1+\sigma)|\yyy | } \prod_{j=1}^N y_j^{\theta \alpha -1}  \Delta(\yyy )^{2\theta }  d\yyy 
\\ \label{:66d}
& =  e^{-|\yyy | } \prod_{k=1}^N y_k^{\theta \alpha -1} \Delta(\yyy )^{2\theta } d\yyy 
.\end{align} 
Equation \eqref{:66c} coincides with \eqref{:36f} for $\zzz=\one $: namely,  $\int \Phi_{\varphi (\yy) }^\theta (\xxx ) \mathsf{m}_{\theta, \alpha }^{N} (d \xxx ) =(\theta(\alpha +N-1))_{\varphi(\yy) }^\theta $.
This also follows from a limiting case of the Kadell's integral formula \cite[Theorem 1 (I)]{Kad97}.
Equation \eqref{:66d} is simply a reformulation of $\zFz (;\yyy, \one )=e^{|\yyy |}$.

\section{Integral transforms for $\theta=1$}\label{s:6}
In this short section, we focus on the case $\theta=1$, and we rewrite our general results in the standard notation given by Schur-type determinantal expressions.
We set $V(\xxx )=\prod_{1\le i<j\le N}(x_i-x_j)=(-1)^{N(N-1)/2}\Delta(\xxx )$, and define the multivariate polynomials
\begin{align*}
 \MMM _{\mm } (\xxx ; a,c )&= \frac{\det [M_{m_j +N-j } (x_i ;a,c)]_{i,j=1}^N  }{V (\xxx ) }
 ,\\
 \LLL _{\mm }^{a} (\xxx  )&= \frac{\det [L_{m_j +N-j }^a (x_i)]_{i,j=1}^N }{V (\xxx )} 
.\end{align*}
Let $s_{\mm } (\xxx )= P_\mm^{1} (\xxx )$ and $s_{\mm }^{\ip } (\xxx )= P_\mm^{1, \ip} (\xxx ) $ denote the Schur and interpolation Schur polynomials, respectively.
We note that 
\begin{align}
 \label{:212a}
s_\mm (\xxx )&= \frac{\det [x_i^{m_j+N-j}]}{V (\xxx )}
,\\
   \label{:212b}
s_\mm (\one )&=\prod_{1\le i<j\le N}^{}\frac{m_i -m_j +j -i}{j-i} = \frac{V (\mm+\delta )} {\delta! } 
,\\ \label{:212c}
s_\mm^{\ip } (\mm +\delta)&= \frac{(\mm+\delta)! }{V (\mm+\delta )}  
.\end{align}

\begin{lemma}\label{l:61}
  For all $\mm\in\PN $ and $\xxx \in \R^N$, we have 
\begin{align*}
M_{\mm}^{1} (\xxx ; a, c)
 = \frac{\delta ! \prod_{j=1}^N (a -N+1)_{j-1}}{(1-c^{-1}) ^{N(N-1)/2}   V(\mm+\delta)  } \MMM _{\mm } (\xxx +\delta ; a-N+1,c )
.\end{align*}
\begin{proof}
From  \cite[Theorem 3.15]{Shi16}, for $\mm,\xx \in \PN$, we have
\begin{align*}
M_{\mm}^{1} (\xx ; a, c)
 &
  =\frac{1}{\delta !} \frac{ \prod_{j=1}^N (a -N+1)_{j-1}}{(1-c^{-1}) ^{N(N-1)/2} s_{\mm}(\one ) s_{\xx}(\one ) } \det[ M_{m_i+N- i } (x_j+N-j ; a-N+1, c)]
.\end{align*}  
Combining this with \eqref{:212b}, we complete the proof.
\end{proof}
\end{lemma}

\begin{lemma}\label{l:62}
For any $\mm\in\PN $,  we have
  \begin{align*}
  \LL_{\mm }^{1, a+N} (\xxx )
   =(-1)^{N(N-1)/2}  V (\mm +\delta ) \LLL _{\mm }^{a}(\xxx )
.  \end{align*}
\begin{proof}
Recall the following identity for $f_j (z) = \sum_{m=0}^\infty a_m(j) z^m $; see \cite[Theorem 1.2.1]{Hua63}:
\begin{align*}
  \frac{\det [f_j (z_i)]}{V (\zzz )}=\sum_{\mm \in \PN } \det [a_{m_j +N-j}(i)] s_{\mm} (\zzz )
.\end{align*}
Using this identity for $f_i(z) :=(1-z)^{-a -1} \exp( -\frac{z}{1-z} x_i ) =\sum_{m\ge 0} L_m^a (x_i)z^m $, we have
\begin{align*}&
    \sum_{\mm \in \PN} \det [ L_{m_j+N-j}^a (x_i) ] s_\mm (\zzz )
  \\&
  =\prod_{k=1}^{N }(1-z_k)^{-a-1}\frac{\det  [\exp( -\frac{z_j}{1-z_j} x_i ) ]}{V (\zzz)}
  \\&
  =\frac{V ( -\frac{\zzz }{1-\zzz })}{V  (\zzz )} \prod_{k=1}^{N} (1-z_k)^{-a-1}  \sum_{\mm \in \PN} \det \Big[ \frac{x_i^{m_j+N-j}}{(m_j+N-j)!} \Big] s_\mm \Big(-\frac{\zzz }{1-\zzz } \Big)
  \\&
  =\frac{(-1)^{N(N-1)/2} V(\xxx )}{\delta !} \prod_{k=1}^{N} (\one -z_k)^{-a-N} {}^{}_{0} F _{0} ^{1} \Big(; \xxx , -\frac{\zzz}{\one -\zzz } \Big)
.\end{align*}
In the last line, we used \eqref{:212b} and \eqref{:212c}.
Hence, by using \eqref{:26c}, we obtain $ \LL_{\mm }^{1, a+N} (\xxx )
   =(-1)^{N(N-1)/2} \delta ! s_{\mm }(\one  ) \LLL _{\mm }^{a}(\xxx ) 
 $. 
Thus, combining this with \eqref{:212b}, we finish the proof.
\end{proof}
\end{lemma}

\begin{proposition}
For $a, \sigma>0$ and $\mm\in\PN$,  we have
\begin{align} \label{:25c}&
 \sum_{\yy \in W_{\ge, d}^N}  \MMM _{\mm } \Big(\yy ; a, \frac{\sigma}{1+\sigma} \Big) \frac{\Delta (\yy )}{\Delta (\sigma \xxx )} \det \Big[\frac{(\sigma x_i)^{y_j} e^{-\sigma x_i } }{y_j !} \Big] 
  =\LLL _{\mm}^{a-1} (\xxx )  \frac{ \sigma^{-\frac{N(N-1)}{2}}   (\mm+\delta)!}{\prod_{j=1}^{N} (a)_{m_j +N-j }}
,  \\ & \notag 
\frac{(1+\sigma)^{aN +N(N-1)/2+|\xx|}}{\prod_{j=1}^{N}\Gamma(a+x_j)}\int_{\nWN }\LLL_{\mm }^{a-1}(\yyy ) \frac{\Delta (\yyy )}{\Delta (\xx )} \det[ y_j^{a+x_i-1} e^{-(1+\sigma) y_j} ] d\yyy 
 \\ &\label{:25d} \qquad \qquad \qquad  \qquad \qquad \qquad 
=\MMM _{\mm} \Big(\xx ;a,\frac{\sigma}{1+\sigma }\Big) \frac{\sigma^{|\mm |+\frac{N(N-1)}{2}}}{(1+\sigma )^{|\mm |}}  \frac{ \prod_{j=1}^{N} (a)_{m_j+N-j}}{(\mm+\delta)!}
.\end{align}
\begin{proof}
Since \eqref{:212a} and \eqref{:212c} yield
\begin{align*}
 \frac{ s_{\varphi (\yy )}(\sigma \xxx ) }{ s_{\varphi (\yy )}^{\ip}(\yy ) } e^{-\sigma |\xxx |}  =  \frac{\Delta (\yy )}{\Delta (\sigma  \xxx )} \det \Big[\frac{(\sigma x_i)^{y_j} e^{-\sigma x_i } }{y_j !} \Big] 
,\end{align*}
we obtain \eqref{:25c} from \eqref{:62a} for $\theta=1$ together with \lref{l:61} and \lref{l:62}.
Similarly, from \eqref{:52a} for $\theta=1$ together with \lref{l:61} and \lref{l:62}, and using \eqref{:212a} and \eqref{:212b},  we have \eqref{:25d}.
\end{proof}
\end{proposition}

\begin{remark}
  Identity \eqref{:25c}, which gives the transform of the multivariate Meixner polynomial by the kernel $\Lambda_{\sigma}^N$, is implicit in \cite[Remark 3.11]{Ass19b}.
\end{remark}

\appendix

\section{Multivariate Charlier and Krawtchouk polynomials}\label{s:A}
In this appendix, we always assume that $\theta >0$.
For a non-negative integer $r$, we extend the shifted factorial to all $x\in\C$ by $(x)_r=\prod_{j=0}^{r-1}(x+j)$, with $(x)_0=1$. 
We introduce multivariate Charlier polynomials $ C_{\mm}^\theta (\xxx ;a)$ for $a\neq 0$ and Krawtchouk polynomials $K _\mm^\theta (\xxx ; p,n)$ for $n\in \Z_{\ge}, p\neq 0$ by
\begin{align}
  C_{\mm}^\theta (\xxx ;a)&:=\sum_{\kk \subseteq \mm }  \frac{1}{ \Psi_{\kk }^{\theta}( \one )}  \binom{\mm }{\kk }_{\theta } \binom{\xxx  }{\kk}_{ \theta } \Big(-\frac{1}{a} \Big)^{|\kk |}   
,  \\
  K _\mm^\theta (\xxx ; p,n)&:=\sum_{\kk \subseteq \mm }  \frac{1}{ (-n)_{\kk}^\theta \Psi_{\kk }^{\theta}( \one  )}  \binom{\mm }{\kk }_{\theta } \binom{\xxx  }{\kk}_{ \theta } \Big(\frac{1}{p} \Big)^{|\kk |}     \quad \text{ for } \mm \subseteq n \one  
.\end{align}
A straightforward computation yields the following relations between the multivariate Meixner, Charlier, and Krawtchouk polynomials.
\begin{proposition}
  We have
  \begin{align}\label{:A1a}
 K _\mm^\theta (\xxx ; p,n)&=M_\mm^\theta \Big(\xxx ; -n, \frac{p}{p-1}\Big) 
, \\ \label{:A1b}
 \limi{s }M_\mm^\theta \Big(\xxx ; s ,\frac{a}{s +a } \Big) &=  C_{\mm}^\theta (\xxx ;a)
,\\ \label{:A1c}
 \limi{s }K_\mm^\theta \Big(\xxx ; \frac{a}{s }, s \Big) &= C_{\mm}^\theta (\xxx ;a)
 . \end{align}
\end{proposition}

\subsection{Basic formulae for multivariate Charlier polynomials}
We list formulae for Charlier polynomials in this subsection.

\begin{theorem}(Double generating functions)\label{p:A4}
  For all $\gamma \in \R $, we have
  \begin{align}
 \label{:A4b}
\sum_{\mm, \xx \in \PN  } (\gamma )_\mm^\theta C_{\mm}^\theta(\xx;a)  \Psi_{\mm}^\theta(a  \uuu ) \Psi_{\xx }^\theta (\zzz ) 
&=\prod_{j=1}^{N}(1-a u_j)^{-\gamma} e^{|\zzz | } \oFz  \Big(
  \begin{matrix} \gamma \\ \, \end{matrix} ; \zzz, \frac{- \uuu }{1-a \uuu } \Big)
,\\
  \label{:A4a}
\sum_{\mm, \xx \in \PN  }   C_{\mm}^\theta(\xx;a)  \Psi_{\mm }^\theta(a \uuu ) \Psi_{\xx }^\theta (\zzz ) 
&
=e^{a | \uuu |} e^{|\zzz |}\zFz (  ;\zzz , -\uuu )
.\end{align}
\begin{proof}
Equation \eqref{:A4b} follows by substituting $\zzz \mapsto \zzz / s$ and $ \uuu \mapsto s \uuu  $ in \eqref{:36b}, and then taking $s\to \infty $ using \eqref{:A1b}.
Similarly, we obtain \eqref{:A4a} by substituting $\zzz \mapsto \zzz / s $ in \eqref{:36a} and taking a limit.
Note that \eqref{:A4b} is derived from \eqref{:27d}, and \eqref{:A4a} is derived from \eqref{:36g}, \eqref{:27c}, or \eqref{:25a}.

\end{proof}
\end{theorem}

\begin{theorem}[Generating functions]\label{p:A3}
 For all $\gamma \in\R$, we have 
\begin{align}\label{A:3a}
\! \sum_{\mm \in \PN }\! (\gamma)_{\mm }^\theta C_{\mm }^\theta(\xx ;a) \Psi_{\mm }^\theta(\uuu ) 
&= 
\prod_{j=1}^{N}(1-u_j)^{-\gamma } \!\! \! \sum_{\mm \in \PN } \! (\gamma )_\mm ^\theta  \binom{\xx }{\mm }_\theta \Big(- \frac{1}{a} \Big)^{| \mm |} \Phi_{\mm }^\theta \Big( \frac{\uuu }{\one -\uuu } \Big)
.\end{align}
In particular, for any $\uuu \in \R^N$, we have
\begin{align}
\label{A:3b}
\sum_{\mm \in \PN } C_{\mm } ^{\theta } (\xx ;a ) \Psi_{\mm }^{\theta } (a \uuu )
&= e^{a |\uuu |} \Phi_{\xx }^\theta(\one  -\uuu )
.\end{align}
\begin{proof}
Equation \eqref{A:3a} follows by substituting $\uuu \mapsto s \uuu $ and using \eqref{:A1b} for \eqref{:34b}.
Equation \eqref{A:3b} follows from \eqref{:34a} or  \eqref{:27a} by substituting $\uuu \mapsto a \uuu / s$ and $s\to \infty$.
\end{proof}
\end{theorem}

Orthogonality, difference equations, and inversion formulae follow directly from \pref{p:39}, \pref{p:313}, and \pref{p:315} as follows.

\begin{proposition}[Difference equations]
For all $\mm\in\PN $ and $\xxx \in \R^N$, we have
  \begin{align*}
-|\mm| C_{\mm}^\theta(\xxx ;a )
&
  =a\sum_{i=1}^N  h_{+,i}^\theta (\xxx )  C_{\mm}^\theta(\xxx +\eee _i ;a)
  \\&\qquad
  -(|\xxx | +aN )C_{\mm}^\theta(\xxx ;a ) +\sum_{i=1}^N h_{-,i}^\theta (\xxx ) (x_i+\theta(N-i))C_{\mm}^\theta(\xxx -\eee _i ;a)
.\end{align*}
\end{proposition}

\begin{proposition}[Orthogonality]
For all $\mm, \kk\in\PN $, we have 
\begin{align*}
  \sum_{\xx \in \PN } e^{-Na}  \Psi_{\xx }^\theta (a \one  ) C_{\mm}^\theta(\xx; a ) C_{\kk}^\theta(\xx; a)
  = \frac{1}{  \Psi_{\mm}^{\theta}(a \one )}\delta_{\mm,\kk }
.\end{align*}
Here, the orthogonality measure is positive if $a>0$.
\end{proposition}

Note that the orthogonalising measure is related to $\Lambda _{\theta, \sigma}^{N}(\xxx, \cdot) $.

\begin{proposition}[Inversion formulae]
  For all $\mm\in\PN$, we have 
\begin{align*}
P_{\mm}^{\theta, \ip } (\xxx +\theta \delta )=  P_{\mm}^\theta (a \one  ) \sum_{\kk \subseteq \mm } (-1)^{|\kk |}\binom{\mm }{\kk }_\theta  C_\kk^\theta (\xxx ;a )
.\end{align*}
\end{proposition}

\subsection{Basic formulae for multivariate Krawtchouk polynomials }
Although many identities for Krawtchouk polynomials arise as limits $a\to -n $ of the corresponding Meixner identities, one has to be careful in this limiting procedure, since the factors $(a)_{\kk}^\theta $ converge to $0$, and may appear in denominators. 
We note that 
\begin{align}
   \! \!
\lim_{a\to -n }(a)_{\mm}^\theta (a)_{\xx}^\theta  M_\mm^\theta \Big(\xx ; a, \frac{p}{p-1} \Big)
 \label{A:8d}
  = \!
  \begin{cases}
    (-n)_{\mm }^\theta (-n)_{\xx }^\theta  K_{\mm}^\theta (\xx ; p,n) &\text{ if }\mm ,\xx \subseteq n \one 
,\\
0 & \text{ otherwise}
,   \end{cases}
\end{align}
since we have 
\begin{align*}
\lim_{a\to -n } \frac{(a)_{\mm} ^\theta } {(a)_{\kk}^\theta } 
=\prod_{j=1}^{N} (-n-\theta(j-1)+k_j )_{m_j -k_j} 
=\frac{(-n)_{\mm } ^\theta }{(-n)_{\kk }^\theta }
\end{align*} 
for $\kk \subseteq \mm $ and $\kk \subseteq n \one$ by a direct computation.
We also note that 
\begin{align*}
  \binom{ \gamma \one }{\kk }_{\theta }=\Psi_{\kk }^\theta  (-\one ) (-\gamma )_{\kk}^\theta
.\end{align*}

\begin{theorem}[Double generating functions]
Assume that $p\neq 1$.
Then we have
\begin{align}\notag &
\sum_{\substack{\mm, \xx \in \PN  \\ \mm ,\xx \subseteq n\one }}  (-n)_\mm^\theta  (-n )_\xx^\theta K_{\mm}^\theta(\xx;p,n)  \Psi_{\mm }^\theta \Big(\frac{p}{p-1}\uuu \Big) \Psi_{\xx }^\theta (\zzz ) 
\\& \notag 
=\sum_{\substack{\mm, \xx \in \PN  \\ \mm ,\xx \subseteq n\one }}  (-1)^{|\mm|+|\xx| } \binom{n\one  }{\mm }_\theta \binom{n \one }{\xx }_\theta K_{\mm}^\theta(\xx;p,n) \Phi_{\mm }^\theta \Big(\frac{p}{p-1}\uuu \Big) \Phi_{\xx }^\theta (\zzz ) 
\\& \label{A:8b}
=\prod_{j=1}^{N}\Big(1- \frac{p}{p-1} u_j \Big)^{n } \oFz \Big(  \begin{matrix} -n \\ \, \end{matrix} ;\zzz , \frac{\one -\uuu }{\one -p(p-1)^{-1} \uuu } \Big) 
.\end{align}
\begin{proof}
Equation \eqref{A:8b} follows by taking $a \to -n $ in \eqref{:36a} with \eqref{A:8d}.
\end{proof}
\end{theorem}

\begin{theorem}[Generating functions]
  For $\xx \subseteq n\one $, we have
\begin{align}\label{A:9b}
\sum_{\substack{\mm  \in \PN \\ \mm \subseteq n \one } } (-n )_{\mm }^\theta K_{\mm } ^{\theta } (\xx ; p , n ) \Psi_{\mm }^{\theta } (\uuu )
&=
\prod_{j=1}^{N}(1-u_j)^{n } \Phi_{\xx }^{\theta }\Big(\frac{1+  \frac{1-p}{p} \uuu }{\one -\uuu } \Big)
. \end{align}
\begin{proof}
  Taking $a\to -n $ in \eqref{:34a}, we obtain \eqref{A:9b}.
\end{proof}
\end{theorem}

Orthogonality, difference equations, and inversion formulae follow directly from \pref{p:39}, \pref{p:313}, and \pref{p:315} by taking $a\to -n $ as follows.

\begin{proposition}[Difference equations]
  For all $\mm \subseteq n\one$ and $\xxx \in \R^N$, we have
\begin{align*}
 |\mm| K_{\mm}^\theta(\xxx ; p,n)
&
  =p \sum_{i=1}^N (x_i -n -\theta(i-1))h_{+,i}^\theta (\xxx ) K_{\mm}^\theta(\xxx +\eee_i ; p,n )
  \\& \qquad 
  -( (2p-1 )|\xxx | -n p N   ) K_{\mm}^\theta(\xxx ;p,n ) 
  \\&\qquad
  +(p-1)\sum_{i=1}^N (x_i+\theta(N-i))h_{-,i}^\theta (\xxx ) K_{\mm}^\theta(\xxx -\eee_i ; p,n)
.\end{align*}
\end{proposition}

\begin{proposition}[Orthogonality]
Assume that  $p\neq 1$.
Then, for all $\mm, \kk\subseteq n\one  $, we have 
\begin{align*}
\sum_{\xx \in \PN } (1-p)^{Nn } (-n)_{\xx}^\theta \Psi_{\xx }^\theta \Big(\frac{p}{p-1} \one  \Big) K_{\mm}^\theta(\xx; p,n ) K_{\kk}^\theta(\xx; p,n ) = \frac{1}{ (-n )_{\mm}^\theta \Psi_{\mm}^{\theta}(p(p-1)^{-1}  \one )}\delta_{\mm,\kk }
.\end{align*}
Here, the orthogonality measure is positive if $0<p<1$.
\end{proposition}

\begin{proposition}[Inversion formulae]
For all $\mm \subseteq n\one  $ and $\xxx \in \R^N$, we have 
\begin{align*}
P_{\mm}^{\theta, \ip } (\xxx +\theta \delta )=(-n)_{\mm}^\theta  P_{\mm}^\theta (-p \one  ) \sum_{\kk \subseteq \mm  } (-1)^{|\kk |}\binom{\mm }{\kk }_\theta  K_\kk^\theta (\xxx ; p,n )
.\end{align*}
\end{proposition}

\section{Estimates for special functions}\label{s:B}
Throughout this section, we assume that $\theta>0$.
Here we give estimates for Jack polynomials and associated functions to justify the convergence of sums and integrals.
\begin{lemma}
    For all $\mm  \in \PN $ and $\xxx \in \R^N $, we have
\begin{align}\label{B:1a}
| \Phi_{\mm }^{\theta}  (\xxx ) |\le  \| \xxx \|^{|\mm |} 
.\end{align}
\begin{proof}
We quote tableau sum formulae for the Jack polynomials \cite[Theorem 6.3]{Sta89}, that is,
\begin{align}\label{B:1c}
  P_{\mm }^{\theta} (\xxx )&=\sum_{T \in \mathrm{SSYT}(\mm ) } \psi_{T} ^\theta \prod_{s=(i,j)\in \mm  } x _{T(s)}  
.\end{align}  
Here, 
\begin{align*}
  \psi_{T}^\theta :=\prod_{ i=1}^{N } \prod_{s \in (R\setminus C) _{m^{(i-1)} / m^{(i)} } } \frac{b_{m^{(i)}}^\theta (s ) }{b_{m^{(i-1)}}^\theta (s ) }
, \qquad 
  b_{\mm }^\theta (s) :=\frac{a_{\mm }(s) +\theta (l_{\mm } (s)+1 )}{a_{\mm }(s) +\theta l_{\mm } (s)+1 }
, \end{align*}
where $a_{\mm }(s) ,l_{\mm } (s)$ are the arm and leg length at $s$, and $T(s)$ denotes the entry in $s$; see also \cite[Section 6.3]{Koo15} for details.
In particular, $\psi_{T}^\theta >0$ for $\theta>0$.
Hence, equation \eqref{B:1a} follows from
\begin{align*}
 | P_{\mm }^{\theta}  (\xxx ) |
 \le \sum_{T  \in \mathrm{SSYT}(\mm  ) } \psi_{T} ^\theta \prod_{s\in \mm }  | x _{T(s)}  |
 \le \sum_{T  \in \mathrm{SSYT}(\mm ) } \psi_{T} ^\theta   \| \xxx \|^{|\mm  |}
 = P_{\mm  }^{\theta}  (\one  ) \| \xxx \|^{|\mm  |}
,\end{align*} 
which completes the proof of (i).
\end{proof}
\end{lemma}

\begin{lemma}
There exists a constant $\Ct{c;B2}:=\cref{c;B2}(\theta, N)>0$ such that, for any $\mm\in\PN$ and $\xxx \in \R^N$, 
\begin{align}\label{B:2b}
 \Psi_\mm ^{\theta} (\xxx ) 
 \le \cref{c;B2} \frac{ \prod_{1\le i< j\le N } (m_i-m_j+\theta (j-i+1) )^{2\theta } }{\prod_{k=1}^N m_k ! } \|\xxx \|^{|\mm |}
.\end{align}
\begin{proof}
From \eqref{:20a} and \eqref{:20b}, we have 
\begin{align*}
 \Psi_\mm ^{\theta} (\one ) 
&
 =
\prod_{k=1}^N \frac{1}{( \theta(N-k)+1 )_{m_k} } \prod_{1\le i<j\le N} \frac{(\theta(j-i+1))_{m_i-m_j}}{(\theta(j-i))_{m_i-m_j}}  \frac{(\theta(j-i)+1)_{m_i-m_j}}{(\theta(j-i-1)+1)_{m_i-m_j}}  
\\&
\le 
\prod_{k=1}^N \frac{1}{m_k! } \prod_{1\le i<j\le N} \frac{\theta(j-i)+m_i-m_j}{\theta(j-i)}  \frac{(\theta(j-i+1))_{m_i-m_j}}{(\theta(j-i-1)+1)_{m_i-m_j}}  
.\end{align*}
Combining this with 
\begin{align*}
\frac{(\theta(j-i+1))_{m_i-m_j}}{(\theta(j-i-1)+1)_{m_i-m_j}}  \le \text{Const.}\times \frac{\Gamma(\theta(j-i-1)+1)}{\Gamma(\theta(j-i+1))} (\theta(j-i+1)+m_i-m_j)^{2\theta -1 }
,\end{align*}
we complete the proof.
Here, we have used the fact that, if $a, b >0$, there exists a constant $c_{a,b} $ such that 
\begin{align*}
 \frac{\Gamma (a +m )}{\Gamma (b +m ) }   
 \le c_{a,b} (a+m )^{a-b} 
\end{align*}
 for $m\ge 0$.
 Indeed, it follows from $\Gamma(a+m)/\Gamma(b+m) \sim m^{a-b}$ as $m\to \infty$.
\end{proof}
\end{lemma}

\begin{lemma}\label{l:B3}
For every compact set $K\subset \R^N$, the following double series converges uniformly in $\zzz \in K$, provided that $\|\uuu\| $ is sufficiently small depending on $K$:
\begin{align}\label{B:10}
  \sum_{\mm \in \PN }\sum_{\kk \subseteq \mm } \bigg| \frac{(\gamma )_{\mm} ^\theta }{(a)_{\mm} ^\theta }  \Psi_{\mm}^{\theta}( \one )   (-1)^{|\kk |} \binom{\mm }{\kk }_{\theta } \frac{(a )_{\mm }^\theta }{(a  )_{\kk} ^\theta} \Phi_{\kk }^{\theta }(\zzz ) \Phi_{\mm }^\theta (\uuu ) \bigg|<\infty
 .\end{align}
  \begin{proof}
Let $c_a=\sup_{\kk\in\PN } |(a)_{\kk }| ^{-1} <\infty$.
    From \eqref{B:1a} and \lref{:31f}, we have
    \begin{align*}
\sum_{\kk \subseteq \mm } \bigg| \frac{1 }{(a)_{\kk } ^\theta }  \binom{\mm }{\kk }_{\theta }  \Phi_{\kk }^{\theta }(\zzz ) \bigg|
 \le 
c_a \sum_{\kk \subseteq \mm }  \binom{\mm }{\kk }_{\theta }  \sup_{\zzz\in K} \|\zzz \|^{\kk}  
=  c_a \Phi_{\mm}^\theta \Big(1+ \sup_{\zzz\in K}\|\zzz \| \Big)
.\end{align*}
Therefore, the left-hand side of \eqref{B:10} is bounded by
\begin{align*}
c_a \sum_{\mm \in \PN } (\gamma )_\mm^\theta \Phi_{\mm }^{\theta }\Big(1+ \sup_{\zzz\in K}\|\zzz \| \Big) \Psi_{\mm }^\theta (\uuu ) 
,\end{align*}
which converges uniformly for sufficiently small $\|\uuu \|$ by \lref{l:B2} (i).
This completes the proof.

  \end{proof}
\end{lemma}

\section*{Declarations}

\noindent
\textbf{Funding} \\
Y. K. is supported by JSPS KAKENHI Grant Numbers 21K13812 and 25K17268.
G. S. is supported by JSPS KAKENHI Grant Numbers 21K13808 and 26K06837.

\bibliographystyle{acm}
\bibliography{reference}

\end{document}